\documentclass[trsc, nonblindrev]{informs3}
\OneAndAHalfSpacedXI 

\usepackage{natbib}
\bibpunct[, ]{(}{)}{,}{a}{}{,}%
\def\bibfont{\small}%
\usepackage{booktabs}

\usepackage{tikz}
\usepackage{pgfplots}

\usepackage{listings}
\usepackage{color}
\usepackage{bold-extra}

\TheoremsNumberedThrough     
\EquationsNumberedThrough    

\definecolor{light-gray}{gray}{0.95}
\usepackage{amsmath}
\usepackage{amsfonts}
\usepackage{mathtools}
\usepackage{mathrsfs}
\usepackage{soul}
\setstcolor{blue}
\usepackage[official]{eurosym}

\usepackage[normalem]{ulem}
\usepackage{graphicx}
\usepackage[singlelinecheck=false, font=scriptsize, labelfont=scriptsize]{caption}
\usepackage{subcaption}
\usepackage{comment}
\usepackage{url}
\usepackage{hyperref}

\usetikzlibrary{positioning}
\usepackage{algorithm}
\usepackage{algpseudocode}
\algnewcommand\algorithmicforeach{\textbf{for each}}
\algdef{S}[FOR]{ForEach}[1]{\algorithmicforeach\ #1\ \algorithmicdo}
\algdef{S}[FOR]{ForAll}[1]{\algorithmicforeall\ #1\ \algorithmicdo}
\newtheorem{Remark}{Remark}
\newcommand\redout{\bgroup\markoverwith
{\textcolor{red}{\rule[.5ex]{2pt}{0.4pt}}}\ULon}
\usepackage{adjustbox}
\usetikzlibrary{decorations.pathreplacing,calligraphy}

\begin{document}
\RUNAUTHOR{Hasturk et al.}
\RUNTITLE{Stochastic Cyclic IRP with Supply Uncertainty}
\TITLE{Stochastic Cyclic Inventory Routing with Supply Uncertainty: A Case in Green-Hydrogen Logistics}
\ARTICLEAUTHORS{%
	\AUTHOR{Umur Hasturk, Evrim Ursavas, Kees Jan Roodbergen}
	\AFF{Department of Operations, Faculty of Economics and Business, University of Groningen, Groningen, the Netherlands\\ \EMAIL{u.hasturk@rug.nl, e.ursavas@rug.nl, k.j.roodbergen@rug.nl}}
	\AUTHOR{Albert H. Schrotenboer}
	\AFF{Operations, Planning, Accounting \& Control Group, School of Industrial Engineering, Eindhoven University of Technology, the Netherlands \EMAIL{a.h.schrotenboer@tue.nl}}
}

\ABSTRACT{Hydrogen can be produced from water, using electricity. The hydrogen can subsequently be kept in inventory in large quantities, unlike the electricity itself. This enables solar and wind energy generation to occur asynchronously from its usage. For this reason, hydrogen is expected to be a key ingredient for reaching a climate-neutral economy. However, the logistics for hydrogen are complex. Inventory policies must be determined for multiple locations in the network, and transportation of hydrogen from the production location to customers must be scheduled. At the same time, production patterns of hydrogen are intermittent, which affects the possibilities to realize the planned transportation and inventory levels. To provide policies for efficient transportation and storage of hydrogen, this paper proposes a parameterized cost function approximation approach to the stochastic cyclic inventory routing problem. Firstly, our approach includes a parameterized mixed integer programming (MIP) model which yields fixed and repetitive schedules for vehicle transportation of hydrogen. Secondly,  buying and selling decisions in case of underproduction or overproduction are optimized further via a Markov decision process (MDP) model, taking into account the uncertainties in production and demand quantities. To jointly optimize the parameterized MIP and the MDP model, our approach includes an algorithm that searches the parameter space by iteratively solving the MIP and MDP models. We conduct computational experiments to validate our model in various problem settings and show that it provides near-optimal solutions. Moreover, we test our approach on an expert-reviewed case study at two hydrogen production locations in the Netherlands. We offer insights for the stakeholders in the region and analyze the impact of various problem elements in these case studies.}
\maketitle

\section{Introduction}

Renewable energy supply is known to be intermittent and uncertain \citep{anvari2016short}. For example, wind and solar energy generation depend on weather conditions \citep{drucke2021climatological}.
A common way to maintain the balance between supply and demand in the electricity network is by means of conventional gas-fired power plants \citep{safari2019natural}, which have an easily adjustable power output. 
The high market prices for natural gas, along with the desire of many countries to eliminate the use of fossil fuels \citep[e.g.,][]{eu2012},
stimulate the development and use of other mechanisms to balance supply and demand of energy in the network. 
The use of hydrogen as a storage medium for energy \citep{liu2010advanced} is such an alternative balancing mechanism, which is expected to be deployed increasingly towards the future \citep{abe2019hydrogen}. 
Hydrogen can be produced emission-free from (renewable) electricity when supply exceeds demand. It can be stored and transported as a gas, and can be converted emission-free into electricity when needed \citep{oliveira2021green}. 
Since conversions between hydrogen and electricity result in significant energy losses, direct use
of hydrogen as a vehicle fuel, a heating source for households, and a feedstock for the industry is also considered attractive \citep{ball2009future}.

A typical hydrogen-based energy network consists of a hydrogen producer and multiple hydrogen customers. 
The producer makes hydrogen from renewable energy sources, which yields a stochastic hydrogen supply per period. 
The producer is responsible for replenishing geographically dispersed customers, who generate electricity from hydrogen, or use it directly.
Customers have stochastic demand for hydrogen per period. 
Although hydrogen can be transported via pipelines, this is for most applications not considered a likely scenario in the near future \citep{staffell2019role}. 
Hydrogen is therefore transported between locations via vehicles.
The producer and the customers each have a storage facility to keep hydrogen in stock to buffer between deliveries. 
The capacities of these storage facilities are limited.
There are two main challenges that arise in the planning for such a network. 

Firstly, we need a transportation delivery schedule for transporting hydrogen from the producer to the customers. Similar to the cyclic inventory routing problem \citep[CIRP, see][]{raa2009practical, rau2018optimization} this consists of: 
(1) a repetitive schedule with fixed time intervals between consecutive deliveries for each of the customers,
(2) the routing for all vehicles, which determines which customers are delivered by which vehicle and in which sequence.
The repetitive schedule and routing are decided upon at the beginning of the planning horizon and are subsequently used in all periods. This helps customers and the producer to align their dependent planning processes and associated resources as they will know exactly when the deliveries will occur.

Secondly, we need to manage hydrogen availability across all locations in response to stochastic supply at the producer and stochastic demand at the customers. As commonly deployed in the stochastic inventory routing problem \citep[SIRP, see][]{raa2021multi, sonntag2021tactical}, the producer determines the delivery quantity for each customer by using an inventory policy. Specifically, for each delivery period --as imposed by the transportation delivery schedule-- we have a base-stock policy for each customer, given a required service level. The base-stock levels for a customer can differ between delivery periods. Furthermore, the limited size of the producers' storage facility, combined with the stochasticity in supply and demand, may  lead to hydrogen shortages or surpluses at the producer's location. The producer is connected to an external hydrogen market at which it can buy and sell hydrogen at market prices to guarantee that customers' service levels are met. Dynamic purchasing decisions at the producer determine when to buy and sell hydrogen, and in which quantity.

It is crucial that both the \textit{static} transportation delivery schedule and the producer's \textit{dynamic} purchasing decisions are optimized jointly to achieve long run global optimal solutions.
This is complex, however, because the static transportation delivery schedule is decided upon only once, while the purchasing decisions are made dynamically in every period that follows. In isolation, these two decision problems are addressed in their own, natural way. Obtaining a static transportation delivery schedule is done using mixed integer programming (MIP) models and methods. On the other hand, the dynamic purchasing decisions are modeled naturally as a Markov decision process (MDP), thus requiring stochastic dynamic programming. In our hydrogen setting, the outcome of the MIP defines the size and timing of the stochastic demand faced by the producer, and thus defines the state, action, and transition space of the MDP. How to best combine this into a joint optimization problem is a yet unresolved fundamental question in inventory routing for which, to the best of the authors' knowledge, no generic solution approach exists.

\begin{figure}[t]
    \centering
    \includegraphics[width=9cm]{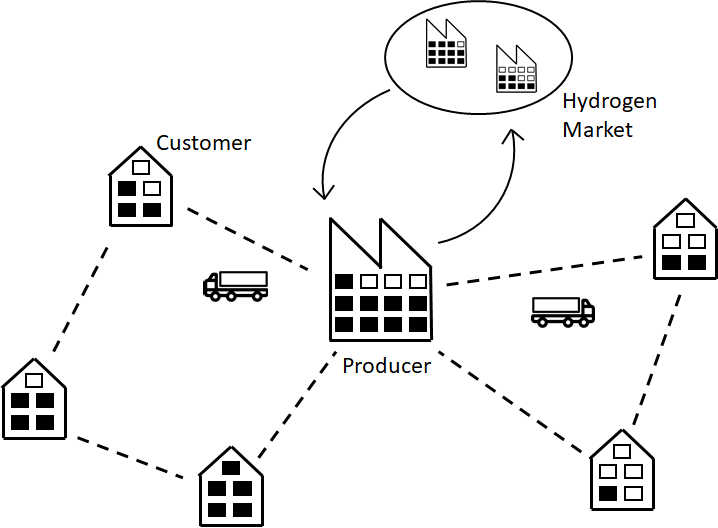}
    \caption{Our paper considers a stochastic, cyclic, inventory routing problem with stochastic supply occurring in hydrogen logistics.} \label{fig:galaxy}
\end{figure}

In this paper, we consider a stochastic  cyclic inventory routing problem (SCIRP) with supply uncertainty, as illustrated in Figure \ref{fig:galaxy}. Its goal is to find a long run, global optimal solution that minimizes the cost of the transportation delivery schedule, the inventory policies for the customers, and the dynamic purchasing decision for the producer. For this, we propose a generic solution approach that iteratively solves MIPs and MDPs utilizing ideas of parameterized cost function approximations \citep{powell2022parametric}. A parameterized cost function approximation is a method for solving stochastic dynamic decision problems by using a modified, parameterized deterministic policy so that it better anticipates future uncertainty. We apply this idea to our static transportation delivery schedule, i.e., we use a modified parameterized MIP to find a transportation delivery schedule that anticipates the impact of the dynamic purchasing decision. We evaluate the dynamic purchasing decision by solving the resulting MDP via value iteration. Our approach includes an efficient search method for setting the parameters in the modified parameterized MIP by iteratively solving the modified MIP and the MDP. To the best of the authors' knowledge, this is the first generic approach that combines such static and dynamic decisions in the context of inventory routing problems. 

This paper makes the following contributions to the literature. First, by introducing and solving the SCIRP with supply uncertainty, we consider a new and emerging application in hydrogen production and distribution that is relevant for an efficient transformation towards a green society. Second, we generalize the existing literature on the SCIRP by considering both a cyclic transportation delivery schedule and stochastic supply. The extant literature typically considers only demand uncertainty in non-cyclic settings \citep{adelman2004price, kleywegt2004dynamic}, supply uncertainty in non-cyclic settings \citep{alvarez2021inventory}, or do not consider supply uncertainty \citep{sonntag2021tactical,  raa2021multi, malicki2021cyclic}. Third, we introduce a generic solution approach based on jointly optimizing dependent static (MIP) and dynamic (MDP) decisions utilizing ideas from parameterized cost function approximations. Our approach finds high-quality solutions in low computation times. Moreover, our approach is generic compared to the extant literature in which convex analysis and enumeration techniques are used to solve or approximate dynamic problem decisions \citep[e.g.,][]{basten2015approximate, mulder2019simultaneous}, our approach does not require such analysis or associated simplifying assumptions. Fourth, we show on a set of new benchmark instances the importance of jointly optimizing static and dynamic decisions, as it can lead up to cost savings of 65\%. Finally, we apply our approach to an extensive case study in the Northern Netherlands, where one of Europe's first Hydrogen Valleys \citep{newenergycoalition2020} is realized, by analyzing how hydrogen should be distributed in expert-reviewed future scenarios.

The remainder of this paper is organized as follows. In Section \ref{sect:lite}, we provide a review of literature that relates to the SCIRP, such as the stochastic IRP, the cyclic IRP, and the production routing problem. In Section \ref{sect:math}, the problem narrative is given, followed by a MIP and an MDP model definition for our joint optimization model. 
In Section \ref{sect:soln}, the proposed solution approach is narrated including a new parameterized MIP model definition. In Section \ref{sect:comp}, the computational experiments are presented. In Section \ref{sect:case}, an expert-reviewed case study demonstrates the impact of various problem elements for a number of future scenarios. In Section \ref{sect:conc}, we conclude our paper. 

\section{Literature Review}  \label{sect:lite}
Our study contributes to the green hydrogen energy literature and several research streams of the inventory routing problem (IRP), which we discuss one by one. First, we discuss the distribution of green hydrogen from an application perspective in Section \ref{literature:Energy}. Next, we discuss the literature on the stochastic IRP and cyclic IRP in Sections \ref{literature:SIRP} and \ref{literature:CIRP}, respectively. We combine those elements into the stochastic cyclic IRP in our paper for which related work is presented in Section \ref{literature:SCIRP}. The literature for combining multiple decision levels of a problem in a joint optimization approach is discussed in Section \ref{literature:joint}. Finally, we discuss how our work differs from other, related, settings in Section \ref{literature:PRP}. 

\subsection{Green Hydrogen Energy Literature} \label{literature:Energy}
Contributions in the green hydrogen literature typically take a high-level supply chain perspective. To the best of the authors'  knowledge, there are no contributions that consider routing and supply uncertainty as in our study. We review the most-related works in the following and refer the interested reader to the review by \cite{sgarbossa2022renewable}. 

\cite{almaraz2014hydrogen} consider a case in France  and minimize environmental and safety objectives while deciding about centralized or decentralized production and storage of hydrogen and different transportation modes. \cite{welder2018spatio} focus on multiple production technologies of green hydrogen, providing a comparison of these technologies for a case study in Germany. None of these works, however, consider inventory routing optimization as in our work.

Instead of considering multiple green-hydrogen production technologies, \cite{woo2016optimization} study the green hydrogen supply chain by focusing solely on biomass production. They are concerned with deciding the allocation of land for yielding multiple crops, which are then used as biomass to produce hydrogen. The authors consider uncertainty in both supply and demand of hydrogen, and touch upon transportation decisions too. They assume a fixed cost of using a transportation mode and of operating it between two locations. Clearly, this does not grasp the importance of routing optimization as the IRP does, and thus does not realistically represent green hydrogen distribution costs as in our study.
 
\subsection{Stochastic Inventory Routing Problem} \label{literature:SIRP}
The stochastic IRP (SIRP) is a collective name for IRPs considering at least one stochastic element. Most often this concerns demand uncertainty, but some recent works also incorporate supply uncertainty and lead time uncertainty. Another important distinction between SIRP models is whether decisions are made dynamically or statically. Dynamic models are built typically upon dynamic programming (DP) and Markov decision process (MDP) techniques, providing real-time operational-level decisions. Static models typically focus on integer programming (IP) techniques, utilizing chance constraints to guarantee some operational performance metric, while focusing on optimizing tactical-level decisions. Our work considers static and dynamic decisions while considering both supply and demand uncertainty. In the following, we only touch upon the most related works, and otherwise refer the interested reader to the excellent review by \cite{coelho2014thirty}.

Some seminal works on dynamic SIRP models exist. \cite{adelman2004price} models the SIRP as an MDP by using approximated value functions with the aim to provide dynamic policies. They formulate two primal-dual formulations, in which the dual feasible solutions are implemented to achieve lower bounds on true value function values and to solve the optimality equations of the MDP model. \citet{kleywegt2002stochastic} also considers a dynamic SIRP and model it as an MDP. The MDP model is provided for the general case where a route may consist of multiple customers. The model is, then, solved by a dynamic programming approximation method for only the direct delivery case, where a single customer is served within a route. Later, \citet{kleywegt2004dynamic} extend this work and solve the general case to near optimality by a decomposition and an approximation of the value function. \cite{bertazzi2013stochastic} study a SIRP with demand uncertainty where each customer has to have a base-stock level equaling their inventory capacity. The authors permit stock-outs with a penalty cost. They model the problem dynamically for a finite time horizon and derive two solution methods. The first is a hybrid rollout algorithm that combines a rollout approach with a deterministic MILP. The second is a branch-and-cut algorithm with additional valid inequalities to solve the problem to optimality. A similar approach is given in \cite{coelho2014heuristics} in which short-term forecasts are injected in a MILP, which is then solved on a rolling horizon basis.

Using integer programming techniques, some authors focus on static models of the SIRP. \cite{yu2012large} model the SIRP with demand uncertainty with chance constraints to achieve a service level while respecting a limited inventory capacity. They present a hybrid solution approach by simplifying the model, linearizing nonlinear elements, and decomposing decisions. They show that a large-scale problem of $200$ customers may be solved to near optimality by their approach. Another SIRP under demand uncertainty with service levels is studied by \cite{crama2018stochastic}, where the product is assumed to be perishable. The authors assume no holding cost unlike most other IRPs, but provide the formulation for the case with a holding cost. They create four solution methods for their formulation of the problem with a finite time horizon and discuss the impact of each element of the problem such as information, storage capacity, and shelf life.

Supply uncertainty has only recently started to gain attention in the SIRP. To the best of our knowledge, \citet{alvarez2021inventory} is the only SIRP paper considering supply uncertainty. The authors modeled the problem as a two-stage integer programming recourse for a finite time horizon with discrete probabilistic supply and demand distributions. In this model, a routing decision is selected at the first step. The optimal transportation delivery schedule is, then, derived for all possible combinations of supply and demand scenarios at the second step. Thus, the delivery decisions may be given after the realization of a scenario, assuming that the decision maker knows the exact amounts of supply and demand for all the periods in the time horizon. Our approach does not rely on generating scenarios and utilizes a chance-constrained model instead. This allows us to guarantee exact service levels for the customers.

\subsection{Cyclic Inventory Routing Problem} \label{literature:CIRP}
The cyclic IRP (CIRP) is the branch of IRPs where the proposed routing solutions are cyclic, and can thus be executed repetitively in an infinite planning horizon. The majority of this research area provides fixed transportation delivery schedules to be implemented repetitively, which helps both customers and producers to plan their (dependent) operations. Most of the literature considers deterministic settings \citep[see, e.g.,][]{michel2012column, diabat2021fixed}.

\cite{gaur2004periodic} study a cyclic IRP to provide a weekly transportation delivery schedule for a well-known supermarket chain in the Netherlands: Albert Heijn. The authors limit the number of customers in each vehicle route by two, which is then improved by a heuristic that allows for more customers in a single route. Their proposed solutions reduced costs by $4\%$ during the first year of implementation. \cite{raa2009practical} work on a CIRP by considering some real-life constraints such as minimal cycle times and daily driving time restrictions. Their column generation-based heuristic includes two main steps; clustering and assigning stock levels, and an insertion heuristic for the routes.

Several studies provide heuristic approaches. \cite{vansteenwegen2014iterated} consider a single vehicle CIRP. The aim is to maximize profits obtained from serving customers. Thus, not all customers need to be delivered. The routing structure is similar to multi-vehicle IRPs since the vehicle may do more than one trip at a period. The authors propose an iterative local search metaheuristic for the problem. Another heuristic approach is discussed in \cite{rau2018optimization} for the CIRP with multiple objectives, minimizing costs and the emission of vehicles. The authors consider both the distance and the weight of the vehicle to calculate emission levels. They propose a variant of the particle swarm optimization algorithm to solve the problem, which reduces the costs and the emission levels by $20\%$ in their instances. In \cite{bertazzi2020exact}, the authors solve a finite-time CIRP with a three-phase approach; combining heuristics with branch-and-cut. The first phase applies a branch-and-cut algorithm. If the optimal solution has not been reached within a preset time limit, a lower bound is taken and the second phase of the heuristic approach is solved so that an upper bound is determined. Branch-and-cut is again applied to solve the problem to optimality in the third phase by using the valid inequalities defined by the first two phases.

\subsection{Stochastic Cyclic Inventory Routing Problem} \label{literature:SCIRP}

The variant of the IRP that relates most closely to our work is the stochastic cyclic IRP (SCIRP), where both customer demand uncertainty is considered while cyclic routing solutions are required. 
\citet{sonntag2021tactical} work on a SCIRP with infinite supply and stochastic demand with known stationary distributions. Their aim is to provide a cyclic replenishment policy with a fixed number of periods in-between consecutive deliveries, where it is guaranteed that each customer is replenished up to their base-stock level. They reformulate the chance-constrained model into an integer programming model by using a Dantzig–Wolfe decomposition. The reformulated model is solved to near optimality for mid- and large-size instances using branch-and-price. Our static (transportation delivery schedule) decisions generalize the approach by \cite{sonntag2021tactical} by adhering to a cyclic pattern of a finite number of days (for example a week). Moreover, they do not consider supply uncertainty.  \cite{raa2021multi} study a  similar problem as \cite{sonntag2021tactical}, and propose a population-based metaheuristic approach, by decomposing the problem into route design and  fleet design stages. Finally, \cite{malicki2021cyclic} address a SCIRP with a finite planning horizon. They allow variable safety stocks and replenishment intervals for the first time in the SCIRP literature. This helps to address non-stationary demand patterns. They model the problem as a mixed-integer chance constraint problem, which is then solved with an adaptive large neighborhood search algorithm that smartly utilizes lot-sizing heuristics.

Our study is the first study that addresses a SCIRP including supply uncertainty. In addition, we detail dynamic purchasing decisions on how to mitigate supply uncertainty using an operational-level MDP to include the associated costs in the design of our cyclic inventory routing design. We propose a generic solution approach for this based on parameterized cost function approximation. To the best of our knowledge, no such generic solution approach exists in the context of vehicle routing or inventory routing.

\subsection{Joint Optimization of Multiple Decision Levels} \label{literature:joint}
Combining decisions on different levels, without the need for using MDPs, is addressed by several authors in various domains. \cite{basten2012joint} and \cite{basten2015approximate} study a complex non-linear optimization problem in the context of maintenance optimization and spare-parts stocking. The authors combine tactical-level maintenance decisions (a so-called Level of Repair Analysis) with an operational-level spare-parts stocking strategy. In \cite{basten2012joint}, the problem is solved to optimality under stylized assumptions and generating trade-off curves. \cite{basten2015approximate} provides a heuristic approach based on a feedback mechanism where solutions to a lower-level optimization problem alter the coefficients of an upper-level optimization problem. Finally, \cite{mulder2019simultaneous} work on maritime shipping by the combination of static and dynamic decisions. A convex analysis is also used in their study to specify the operational-level decisions.

The genericity of our approach stands out when comparing it to this stream of literature. Our approach is similar in the sense that we also rely on iteratively solving two models. 
These existing works, however, optimize operational-level decisions within the context of a set of predefined rules and policies (i.e., smart enumeration). Our generic approach evaluates and solves Markov decision processes in each iteration to obtain state-dependent optimal decisions. 
It does not require for any assumptions on the relationship between the MDP and the set-partitioning formulation since it obtains a combined decision both for tactical-level static (MIP) and for operational-level dynamic (MDP) levels.

\subsection{Other Related Research} \label{literature:PRP}

The production routing problem (PRP) concerns jointly optimizing lot sizing and vehicle routing. Unlike IRPs with stochastic supply (as in our problem), PRPs plan the supply (i.e., the lot sizing) based on the routing policy so that supply, inventory holding, and transportation costs are minimized. Studies are mostly based on heuristic multi-phase approaches, but there are a few exact methods available. We only review the most important contributions in this area due to PRP's focus on discrete scheduling of production processes, and refer the interested reader to \cite{adulyasak2015production} and \cite{neves2019solving} for an overview of the complete field.

Regarding the multi-phase approaches, \cite{adulyasak2014optimization} work on a deterministic PRP under a finite time horizon and propose an adaptive large neighborhood search heuristic. This decomposition-based heuristic has two phases; initialization and improvement, where a network flow subproblem is solved in the improvement phase. \cite{adulyasak2015benders} use a similar two-phase construction to solve a stochastic PRP with demand uncertainty. Benders decomposition is used to solve this problem to optimality, enhanced with additional valid inequalities. \cite{neves2019solving} consider a multi-product PRP with time delivery windows. The authors' heuristic approach adds one more phase, thus having three phases named as; size reduction, initial solution, and improvement. This decomposition method is tested in a case study in a meat industry company, achieving a $20\%$ cost reduction compared to the company's policy.

\cite{schenekemberg2021two} introduce the two echelon PRP. The authors provide two exact solution methods; branch-and-cut and a parallel algorithm. The parallel algorithm provides previously unknown optimal solutions for large-sized PRP instances with $50$ customers, and is also found to be efficient for the deterministic two-echelon multi-depot inventory routing problem.

\section{Problem Definition} \label{sect:math}
Here we present the mathematical models for our stochastic cyclic inventory routing problem (SCIRP) with supply uncertainty. We discuss three models; one for the tactical-level inventory decision, one for the operational-level purchasing decision, and one for the joint optimization of both levels. We start with the problem narrative in Section \ref{sect:narr}, where on a high-level the system elements, the decisions to be taken, and the associated costs are explained. Afterward, we introduce a MIP for the static, tactical-level inventory routing decisions in Section \ref{sect:modelIRP} and an MDP for the dynamic, operational-level purchasing decisions in Section \ref{sect:modelMDP}. The joint optimization problem that formalizes the goal of the SCIRP is presented in Section \ref{sect:modelJoint}. In Section \ref{sect:soln}, we provide our cost function approximation approach including a parameterized MIP to solve the problem presented in this section efficiently. Although our focus is on the distribution of green hydrogen, our underlying modeling concept has wider applicability to other inventory routing settings with stochastic supply. This includes settings where inventory availability at the warehouse is stochastic due to complex dependent production planning processes or where inventory is shipped from external suppliers that are subject to disruptions due to labor shortages or stochastic travel times. Nevertheless, each of such examples would require small changes in the modeling. Therefore, we present the model in terms of green hydrogen distribution to keep our exposition focused.

\subsection{Problem Narrative} \label{sect:narr}
We consider a joint green hydrogen production and distribution system over an infinite discrete time horizon that consists of cyclic repetitions of $T$ periods (e.g., a week), where $T$ is given. Every period within the infinite time horizon corresponds to a cycle period $t \in \mathscr{T} \coloneqq \{1, \ldots T\}$. 
A green hydrogen producer receives a stationary stochastic supply from a renewable energy source at each period. 
The producer uses the stochastic supply to replenish a given set of geographically dispersed customers that each face stationary 
stochastic demand at each period. 
We define the node set $\mathscr{N}^0 \coloneqq \{0\} \cup \mathscr{N}$, where $0$ represents the producer and $\mathscr{N} \coloneqq \{ 1,2,\ldots,N \}$ represents the customers.

The decisions comprise two dependent parts. Firstly, a tactical-level inventory routing decision is made in the form of a periodic, routing solution that repeats itself every $T$ periods over the infinite time horizon. These routes are used for the replenishment of hydrogen from the producer to the customers using capacitated vehicles. To design such a so-called transportation delivery schedule, customers are grouped into fixed, mutually exclusive and collectively exhaustive clusters with an associated vehicle route and a delivery schedule over the $T$ periods. Thus, each customer is part of exactly one cluster. The delivery schedule encodes at which of the $T$ periods the associated vehicle route is performed. This  unique vehicle route for a cluster is used for all deliveries of the cluster over the cycle. Customers are delivered at most once every period.
We assume that customers' inventory is replenished according to a base-stock policy subject to a service level constraint at each delivery period. The customer demand is back-ordered  until the next delivery period in a case of stock-out at the customer's location. Note, the base-stock level set for each delivery period is part of our decision. As a consequence of the base-stock levels set and stochastic demand, the amount of inventory on each vehicle route is stochastic. The proposed clustering of customers and associated transportation delivery schedule transforms individual stochastic customer demands into $T$ independent demand distributions for each cycle period for the producer. See Section \ref{sect:modelIRP} for the details. 

Secondly, the producer faces a dynamic, operational-level purchasing problem to manage, at each period, the inflow of stochastic hydrogen supply to ensure the feasibility of the outflow of demand as depicted by the transportation delivery schedule outlined above. 
The producer observes how much demand needs to be replenished by the customers at each period. The producer has access to a capacitated storage unit to keep hydrogen in stock to anticipate the fulfillment of demand in later periods. Furthermore, the producer is connected to a hydrogen market at which it can buy hydrogen to meet customer demand, and sell hydrogen for an extra profit. We ensure in this way that the producer always replenishes inventory up to the decided base-stock policies (as part of the tactical-level decision) at the customers. See Section \ref{sect:modelMDP} for the details.

The goal of the SCIRP is to minimize the total expected cycle costs of the combined tactical-level and operational-level decisions. This includes transportation costs (fixed and variable), inventory holding costs at customers, emergency shipment costs (in case of insufficient vehicle capacity), and hydrogen buying and selling costs. The tactical-level MIP and the operational-level MDP are discussed in Sections \ref{sect:modelIRP} and \ref{sect:modelMDP}, respectively. The ultimate framework of our solution approach is given in Figure \ref{fig:approachFrame}. Note that this involves a cost-function approximation (see Section \ref{sect:soln}) based upon the problem-defining MIP that we introduce below.

\subsection{The Tactical-Level Inventory Routing} \label{sect:modelIRP}
We propose a set-partitioning model for the tactical-level inventory routing decisions. We consider a set of clusters $\mathscr{R}$, where each cluster $r \in \mathscr{R}$ describes a set of customers $\mathscr{N}_r \subseteq \mathscr{N}$, the shortest vehicle route among these customers, the set of delivery periods $\mathscr{T}_r \subseteq \mathscr{T}$, and a base-stock level at each delivery period for each customer. We refer to these delivery periods and the associated base-stock levels as a \textit{transportation delivery schedule} of a cluster. The tactical-level inventory routing decision, then, boils down to selecting clusters with mutually exclusive customers that jointly partition the customer set. An overview of notation introduced in this section is given in Table \ref{tab:NOTATION1}.

\begin{table}[h]
\centering
\caption{Overview of notation for the tactical level} \label{tab:NOTATION1}

\begin{tabular}[h]{p{0.08\linewidth}  p{0.88\linewidth}} 
\toprule
\multicolumn{2}{l}{Sets:}                                                                                                                                                    \\
$\mathscr{T}$                   & set of cycle periods                                                                                                                           \\
$\mathscr{T}_r$ & set of delivery periods of cluster $r$ \\
$\mathscr{N}$                   & set of customers                                                                                                                           \\
$\mathscr{N}^0$                 & set including all customers and the producer                                                                                               \\
$\mathscr{N}_r$                 & set of customers in cluster $r$                                                                                                            \\

$\mathscr{R}$                   & set of clusters                                                                                                                            \\
$\mathscr{\overline{R}}$ & set of selected clusters in a solution \\\\
\multicolumn{2}{l}{Indices:}                                                                                                                                                 \\
$t$, $\ell$                             & indices for cycle periods                                                                                                                        \\
$T$                             & number of cycle periods                                                                                                                        \\
$i$                             & index for nodes                                                                                                                            \\
$N$                             & number of customers                                                                                                                        \\
$r$                             & index for clusters                                                                                                                         \\                                                                                                                             \\
\multicolumn{2}{l}{Parameters:}                                                                                                                                              \\
$\beta_r^i$                     & parameter that equals $1$ if cluster $r$ contains customer $i$, $0$   otherwise                                                            \\
$d_r$                           & the vehicle route length of cluster $r$                                                                                                    \\
$IS_{r}^{it}$      & random variable of inventory position on customer $i$ of cluster $r$ at start of period $t$              
\\
$IE_{r}^{it}$      & random variable of inventory position on customer $i$ of cluster $r$ at end of period $t$                                        
\\
$q^{it}_r$ & random variable of shipped units to customer $i$ of cluster $r$ at period $t$ 
\\
$Q$                             & vehicle capacity                                                                                                                           \\
$U$                             & customer inventory capacity                                                                                                                \\
$\alpha$                        & target service-level                                                                                                                       \\
$\gamma$                        & target probability that the vehicle capacity is not exceeded                                                                               \\
$\mu_i$                         & mean of the supply/demand distribution of node $i$                                                                                         \\
$\sigma_i$                      & standard deviation of the supply/demand distribution of node $i$                                                                           \\
$F_{in}$                        & cumulative distribution function of the demand faced by customer $i$ during $n$ periods 
\\
$D_{in}$ & random variable of the cumulative demand faced by customer $i$ for $n$ consecutive periods \\
$z_\alpha$ & z-score of $\alpha$ \\
$n_{tr}$ & the number of periods until the next replenishment for a delivery period $t$ of cluster $r$ \\
$m_{tr}$ & the number of periods since the previous replenishment for a delivery period $t$ of cluster $r$ \\                                                              \\
\multicolumn{2}{l}{Decision variables:}                                                                                                                                    \\
$x_r$                           & binary decision variable that equals $1$ if cluster $r$ is in the solution, $0$ otherwise                                                \\

   \\
\multicolumn{2}{l}{Cost   components:}                                                                                                                                       \\
$c^{ \textsc{t}}_r$             & cyclic transportation cost of cluster $r$                                                                                                  \\
$c^{\textsc{h}}_r$              & expected cyclic holding cost of cluster $r$                                                                                                         \\
$c^{\textsc{e}}_r$              & expected cyclic emergency shipment cost of cluster $r$                                                                                              \\
$c_r$                           & expected cyclic tactical cost of cluster $r$                                                                                         \\
$W$                             & fixed transportation cost of a replenishment                                                                                               \\
$w$                             & variable transportation cost per distance                                                                                                  \\
$h$                             & holding cost per unit per period                                                                                                           \\
$e$                             & emergency shipment cost per unit                               \\ 

\bottomrule     
\end{tabular}

\end{table}

We assume a base-stock policy at each delivery period $t \in \mathscr{T}_r$ for each cluster $r \in \mathscr{R}$ with associated transportation delivery schedule, covering the demand with service level $\alpha$ until the next delivery period for the customers in $\mathscr{N}_r$. 
For a delivery period $t$ and a cluster $r$, let $n_{tr}$ be the number of periods until the next delivery and let $m_{tr}$ be the number of periods since the last delivery. We assume the replenishment occurs instantaneously at the start of the period. Let $F_{in}$ denote the cumulative distribution function of the demand faced by customer $i \in \mathscr{N}$ during $n$ periods. We assume demand (and later supply) are continuous random variables. We set the base-stock level equal to $F^{-1}_{i(n_{tr})}(\alpha)$ to achieve a service level $\alpha$ at each customer level, i.e., the critical fraction of demand distribution over the next $n_{tr}$ periods at customer $i \in \mathscr{N}_r$ at delivery period $t$. For example, if demand is normally distributed with mean $\mu_i$ and standard deviation $\sigma_i$,  $F^{-1}_{i(n_{tr})}(\alpha)$ equals $n_{tr} \mu_i + z_\alpha \sigma_i \sqrt{n_{tr}}$ where $z_\alpha$ is the z-score associated with the $(1-\alpha)$ quantile.

To properly define the inventory dynamics at each customer, we introduce the random variable ${IS}_r^{it}$ representing the inventory position at customer $i \in \mathscr{N}_r$ at the start of each period $t$ of cluster $r \in \mathscr{R}$. Similarly, let the random variable ${IE}_r^{it}$ representing the inventory position at customer $i \in \mathscr{N}_r$ at the end of each period $t$ of cluster $r \in \mathscr{R}$. A holding cost of $h$ is paid per unit per period, where the unit is calculated as the average of the start and the end inventory positions of a period.

Let $q_r^{it} \coloneqq F^{-1}_{i(n_{tr})}(\alpha) - IE_r^{i,t-1}$ be the stochastic quantity to be shipped to customer $i$ as part of cluster $r$ on delivery day $t$, and $q_r^{it} \coloneqq 0$ for non-delivery days. For the example normal distribution case, $IE_r^{i,t-1} \sim N( z_\alpha \sigma_i \sqrt{m_{tr}} , m_{tr} \sigma_i^2 )$ and thus $q_r^{it} \sim N( n_{tr} \mu_i + z_\alpha \sigma_i ( \sqrt{n_{tr}} - \sqrt{m_{tr}} ), m_{tr} \sigma_i^2) $ on delivery day $t$. Then, the stochastic total quantity $\sum_{i \in \mathscr{N}_r}q_r^{it}$ to be replenished on delivery day $t$ by cluster $r$ potentially exceeds the vehicle capacity $Q$. An emergency shipment is called in this case to deliver the remaining amount at a unit penalty cost $e$. Finally, the fixed transportation cost is $W$ per replenishment, and the variable transportation costs equal $w \cdot d_r$, where $d_r$ is the length of the shortest route among the customers in cluster $r$.

The cyclic cluster cost $c_r$ for each $r \in \mathscr{R}$ consists of the deterministic fixed and variable transportation costs ($c^{ \textsc{t}}_r$), the expected inventory holding costs at its customers ($c^{\textsc{h}}_r$), and the expected costs of emergency shipments ($c^{\textsc{e}}_r$). The cluster cost $c_r$ is then defined by
\begin{align}
     c_r & = c^{\textsc{t}}_r + c^{\textsc{h}}_r + c^{\textsc{e}}_r \\ 
    & = \sum_{t \in \mathscr{T}_r}  \left[ \left(W + w\cdot d_r\right) +  h \sum_{i \in \mathscr{N}_r} \sum_{\ell = t}^{t + n_{tr}-1} \mathbb{E}\left[ \frac{1}{2}( (IS_{r}^{i\ell})^+ + (IE_r^{i\ell})^+) \right] + e \cdot \mathbb{E} \left[ \left( \sum_{i \in \mathscr{N}_r}  q^{it}_r - Q\right)^+ \right] \right],   \label{clusterCostFunction}
\end{align}
where $(\cdot)^+$ denotes $\max(0,\cdot)$. The first term of Equation \eqref{clusterCostFunction} represents the cycle transportation cost for a cluster consisting of fixed and variable cost. The second term represents the expected cycle holding costs for a cluster by taking expectation of the average inventory level, and the third term denotes the expected cyclic emergency shipment costs. Emergency shipments involve only small quantities, and aligning with the insights provided by our industry practitioners, they only include a variable (though higher) cost.

Each cluster $r \in \mathscr{R}$ is subject to three constraints. First, each customer has only a limited capacity $U$ for the storage of hydrogen, so $F^{-1}_{i(n_{tr})}(\alpha) \leq U$ should be satisfied for all $r \in \mathscr{R}, i \in \mathscr{N}_r, t \in \mathscr{T}_r$. Second, by setting base-stock levels equal to $F^{-1}_{i(n_{tr})}(\alpha)$, we implicitly impose that the probability of having no back-orders is at least $\alpha$ for each customer $i$ in cluster $r$ on period $t$. 
\begin{align}
\mathbb{P} \Big\{ IE_{r}^{it} \geq  0 \Big\}  \geq \alpha \quad\quad \forall r \in \mathscr{R}, i \in \mathscr{N}_r, t \in \mathscr{T}. \label{cc1}
\end{align}

Furthermore, Constraint \eqref{cc2} limits the probability that vehicle capacity is not sufficient (i.e., the probability an emergency shipment occurs) with at most $1-\gamma$ for each delivery of each cluster $r$. The chance constraint helps to limit the size of the cluster set $\mathscr{R}$, which is computationally convenient.
\begin{align} 
\mathbb{P} \Bigg\{ \sum_{i \in \mathscr{N}_r} q_r^{it} \leq Q \Bigg\} \geq \gamma \span \quad\quad \forall r \in \mathscr{R}, t \in \mathscr{T}_r. \label{cc2} 
\end{align}

By setting the base-stock level as $F^{-1}_{i(n_{tr})}(\alpha)$, we act as if we have a classic single-item inventory system with back-orders and lead time equal to the number of periods that we need to cover till the next replenishment, and we act as if this repeats continuously. This is optimal in case the number of periods between deliveries is equal (see also Lemmas 1-3 in \cite{sonntag2021tactical}). However, in practice, this is often not possible. For example, if the cycle $\{1, \ldots, T\}$ is a week and two delivery days are selected, the number of days between consecutive deliveries is not equal. Thus, if $n_{tr}$ and $m_{tr}$ differ, we obtain a slightly higher service level than $\alpha$. 
The following remarks make this explicit. 

\begin{Remark}[Order-size distribution] \label{remark1}
We approximate $q_r^{it} \coloneqq \left( F^{-1}_{i(n_{tr})}(\alpha) - IE_r^{i,t-1} \right)^+$ as $q_r^{it} \coloneqq F^{-1}_{i(n_{tr})}(\alpha) - IE_r^{i,t-1}$ in this paper. Let $t \in \{1, \ldots, T\}$ be a delivery period of customer $i$, assume $m_{tr} >> n_{tr}$ and $T$ is sufficiently large, and assume only 2 delivery periods are selected for customer $i$ within the cycle $\{1, \ldots, T\}$. In such a case, it is non-negligible that the observed inventory level on period $t$ exceeds the set base-stock policy 
(i.e., $IE_r^{i,t-1} > F^{-1}_{i(n_{tr})}(\alpha) $), and accordingly, no order will be placed at the producer. That is, the order size distribution at time $t$ will have a significant probability mass at $0$. 
\end{Remark}

\begin{Remark}[Optimality of base-stock policy]
Remark \ref{remark1} deduces that setting base-stock levels as we do might result in overshooting the set base-stock levels. If this is the case, it will result in a too-high service level (i.e., exceeding $\alpha$) at a customer. This follows trivially by realizing that the total safety stock increases if we replenish more often (i.e., demand pooling) during the cycle. Thus, Constraint \eqref{cc1} is respected but not necessarily binding in our experiments. However, preliminary experiments have shown that these effects were not significant from a computational perspective. That is, it will not lead to different routing decisions. We, therefore, ignore this phenomenon in the remainder of this paper.
\end{Remark}

\begin{Remark}[On the emergency shipment cost]
The approximation on $q_r^{it}$ overshoots emergency shipment costs. But while approximation error is seen on cases when $m_{tr} >> n_{tr}$, emergency shipment costs are significant when $m_{tr} << n_{tr}$. As also the preliminary experiments have shown, this phenomenon almost never happens, and thus ignored in the remainder of the paper.
\end{Remark}

In what follows, we detail the calculation of holding costs. In line with the previous remarks, these expressions are thus approximations of the actual cost, but preliminary experiments also show that the approximations are good enough.

Recall that the holding cost is calculated by the average inventory over the cycle periods, assuming the demand is observed uniformly throughout a period. Let $D_{in}$ be the random variable for the cumulative demand faced by customer $i$ during $n$ periods. On $(t + \ell)^{th}$ period $(\ell\leq n_{tr})$, after the last delivery on cycle period $t$, the average inventory on customer $i$ with the base-stock policy is 
\begin{align}
     \frac{1}{2} \left( \mathbb{E}\left[ \left( F^{-1}_{i(n_{tr})}(\alpha) - D_{i\ell} \right)^+ \right] + \mathbb{E}\left[ \left( F^{-1}_{i(n_{tr})}(\alpha) - D_{i,\ell+1} \right)^+ \right] \right).  \label{holdNormal}
\end{align}
For the example of normally distributed demand, we have $ F^{-1}_{i(n_{tr})}(\alpha) - D_{i\ell}  \sim N\left( \left(n_{tr} -\ell \right) \mu_i\right)+ z_\alpha \sigma_i \sqrt{n_{tr}}, \ell \sigma_i^2)$. From this, we can calculate the expected holding cost for cluster $r$ as following:
\begin{align}
    c^{\textsc{h}}_r & =  h \sum_{t \in \mathscr{T}_r} \sum_{i \in \mathscr{N}_r} \sum_{\ell = t}^{t + n_{tr}-1} \mathbb{E}\left[ \frac{1}{2}( (IS_{r}^{i\ell})^+ + (IE_r^{i\ell})^+) \right] \\
    &= h \sum_{t \in \mathscr{T}_r} \sum_{i \in \mathscr{N}_r} \sum_{\ell = 0}^{n_{tr}-1}  \frac{1}{2} \left( \mathbb{E}\left[ \left( F^{-1}_{i(n_{tr})}(\alpha) - D_{i\ell} \right)^+ \right] + \mathbb{E}\left[ \left( F^{-1}_{i(n_{tr})}(\alpha) - D_{i,\ell+1} \right)^+ \right] \right).\label{eq:holding}
\end{align}
In Equation \eqref{eq:holding}, we calculate the average inventory per period between the replenishments and penalize that with the holding cost. The expectations can be made explicit by standard techniques and thus \eqref{eq:holding} can be evaluated exact.

The tactical-level inventory routing decision is, then, obtained by solving a set-partitioning model. Let $\beta_r^i$ be equal to $1$ if  cluster $r$ contains customer $i$, and 0 otherwise. Let $x_r$ be a binary decision variable equaling 1 if cluster $r$ is selected and 0 otherwise. The tactical-level inventory routing decision is then obtained by solving:
\begin{align}
    \text{min  } & \sum_{r \in \mathscr{R}} c_r x_r\label{initStart}  \\
    \text{s.t  } & \sum_{r \in \mathscr{R}} \beta_r^i x_r = 1 & \forall i \in \mathscr{N}, \label{custAssignBase}\\
    & x_r \in \{0,1\} & \forall r \in \mathscr{R}. \label{initEnd}
\end{align}
The objective is to minimize the expected cyclic costs of selected clusters. Constraints \eqref{custAssignBase} ensure each customer is contained in exactly one cluster. This model relies on a full enumeration of the cluster set $\mathscr{R}$. We detail how we do this in Section \ref{sect:soln}.
The solution to the tactical-level inventory routing problem comprises a set of selected customer clusters, which we will denote by $\mathscr{\overline{R}} \coloneqq \{ r \in \mathscr{R} \mid x_r = 1\}$. As the tactical-level inventory routing model ignores supply uncertainty faced by the producer, the next section will introduce a dynamic purchasing model for the producer that takes as input the selected set of clusters $\mathscr{\overline{R}}$ to ensure the feasibility of tactical-level inventory routing decisions in operational settings. Finally, we remark that our MIP model can easily be extended by other side constraints, for example on the number of clusters per day. This might help to reduce the impact of supply uncertainty for some particular parameter settings, but this will remain rather dependent on the parameter values and instance characteristics and reduces the genericity of our model and approach. Other appropriate side constraints, however, can easily be added to the formulation if needed in some given real context.

\subsection{The Operational-Level Purchasing} \label{sect:modelMDP}

The producer needs to ensure at each period that enough hydrogen is available to satisfy customer demand as imposed by the transportation delivery schedules of the selected clusters. The transportation delivery schedule is the result of the tactical-level inventory routing decision, as specified in Section 3.2. The producer has a maximum capacity to store hydrogen, and it can buy and sell hydrogen at fixed market prices at the hydrogen market. 
We model this dynamic purchasing decision of how much hydrogen to buy and sell at each period, and thus how much hydrogen to keep in stock at the producer, as an MDP.  
We first define the order of events at each period, after which we provide the MDP formulation. An overview of notation introduced in this section is given in Table \ref{tab:NOTATION2}.

\begin{table}[h]
\centering
\caption{Overview of notation for the operational level} \label{tab:NOTATION2} 

\begin{tabular}[h]{p{0.15\linewidth}  p{0.81\linewidth}} 
\toprule
\multicolumn{2}{l}{Sets:}   \\ 
$\mathscr{A}(s_2)$ & set of admissible buying and selling decisions on state $s_2$\\
$\Phi_t$ & set of net outflow levels of period $t$
\\      
$\Pi(\mathbf{x})$ & set of feasible dynamic purchasing policies on tactical-level solution $\textbf{x}$ \\\\
\multicolumn{2}{l}{Indices:}      \\                     
$\omega_1$, $\omega_2$ & indices for state inventory positions \\
$s_1,s_2$ & indices for states
\\                                                          $\phi_t$ & index for net outflow levels of period $t$ \\                                                                        \\
\multicolumn{2}{l}{Parameters:}                                                                                                                                        \\
$O_t$ & random variable of net outflow quantity of period $t$\\ 
$O^+_t$ & random variable of replenishment quantity of period $t$\\ 
$O^-_t$ & random variable of supply quantity of period $t$\\ 
$\overline{U}$ & producer inventory capacity\\\\
\multicolumn{2}{l}{Decision variables:}                                                                                                                                    \\
$q_1$, $q_2$ & amount of purchased and sold inventory\\
$\mathbf{x}$ & tactical-level inventory routing solution\\
$\pi(\mathbf{x})$  & the policy of dynamic purchasing problem, depending on tactical solution $\mathbf{x}$\\
   \\
\multicolumn{2}{l}{Cost   components:}                                                                                                                                       \\ 
$c(\omega_2,q_1,q_2)$ & cost of decision $(q_1,q_2)$ with inventory position of $\omega_2$\\
$v(s_1), y(s_2)$ & future costs of states \\ 
$b_1,b_2$ & variable buy/sell costs of the hydrogen \\
$K_1,K_2$ & fixed costs of purchases from the outside market \\
$C(\pi(\mathbf{x}))$ & expected cyclic purchasing cost of the policy $\pi(\mathbf{x})$\\ 

\bottomrule     
\end{tabular}

\end{table}

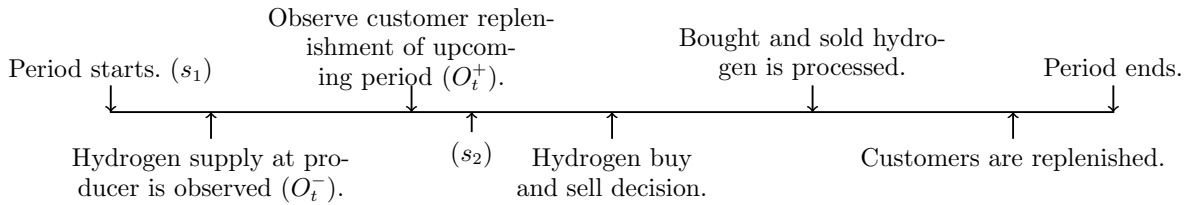
\begin{figure}[h]
\raggedright
\resizebox{18cm}{!}{
\begin{tikzpicture}[x=3cm,nodes={text width=5cm,align=left}]
\draw[black,-,thick,>=latex,line cap=rect]
  (-0.5,-0.3) -- (4.5,-0.3);
\foreach \Xc in {-0.5, 1, 3, 4.5}
{
  \draw[black,thick, ->] 
    (\Xc,0.1) -- (\Xc, -0.3);
}
\foreach \Xc in {0,2,4}
{
  \draw[black,thick, ->] 
    (\Xc,-0.7) -- (\Xc, -0.3);
}
\node[below,align=center,anchor=north,inner xsep=0pt,color=black] 
  at (-0.5,0.65) 
  {Period starts. ($s_1$)};  

\node[below,align=center,anchor=north,inner xsep=0pt,color=black] 
  at (0,-0.7) 
  {Hydrogen supply at producer is observed ($O^-_t$).};  

\node[below,align=center,anchor=north,inner xsep=0pt] 
  at (1,1.4) 
  {Observe customer replenishment of upcoming period ($O^+_t$).};  

\node[below,align=center,anchor=north,inner xsep=0pt] 
  at (2,-0.7) 
  {Hydrogen buy and sell decision.};

\node[below,align=center,anchor=north,inner xsep=0pt] 
  at (3,1.1) 
  {Bought and sold hydrogen is processed.};

\node[below,align=center,anchor=north,inner xsep=0pt] 
  at (4,-0.7) 
  {Customers are replenished.};

  \node[below,align=center,anchor=north,inner xsep=0pt,color=black] 
  at (4.5,0.65) 
  {Period ends.};  


  \node[below,align=center,anchor=north,inner xsep=0pt,color=black] 
  at (1.3, -0.6) 
  {$(s_2)$}; 
  \draw[black,thick, ->] 
    (1.3,-0.6) -- (1.3, -0.3);
  
\end{tikzpicture}

}
\caption{Timeline of a period.}
\label{fig:timeline}
\end{figure}

At each period, we first observe the current hydrogen inventory level and the customer replenishment within the clusters that are delivered in that period (as given by the tactical-level decision). Hereafter, the producer buys or sells hydrogen from the market, and we assume this arrives instantaneously (e.g., overnight) at the producer. Finally, the customers are replenished. A graphical overview is presented in Figure \ref{fig:timeline}.

The clusters $\mathscr{\overline{R}} \coloneqq \{r \in \mathscr{R} \mid x_r = 1\}$, determined by the tactical-level inventory routing problem, partly define the MDP formulation. These clusters impose at each period $t$ a total stochastic replenishment quantity $O^+_t$ that the producer needs to satisfy due to the base-stock policies set at the customers associated with the clusters $\mathscr{\overline{R}}$, where $O^+_t = \sum_{r \in \mathscr{\overline{R}}} \sum_{i \in \mathscr{N}_r} q_r^{it}$. 
In addition, the producer faces stochastic supply $O^-_t$ of green hydrogen from renewable energy sources at each period $t$. We assume that the distribution of $O^-_t$ is known. Let $O_t$ be the net stochastic outflow quantity, i.e., $O_t := O^+_t - O^-_t$.

We consider two states in each period in the MDP. First, the state $s_1$ is observed before the realization of the net stochastic outflow quantity $O_t$, and consists of the period $t \in \{1, \ldots, T\}$ and the inventory level $\omega_1 \in \mathbb{Z}_{\geq 0}$. The state $s_2$ is observed after the realization of $O_t$, and consists of the period $t$ and the inventory position $\omega_2 \in \mathbb{Z}$. The state $s_1 = (t, \omega_1)$ transitions to the state $s_2 = (t, \omega_2)$, where $\omega_2 = \omega_1 - O_t$. 
In state $s_2 = (t, \omega_2)$, we transition to state $s_1 = (t + 1 \mod T,\omega_1)$ by taking hydrogen buy ($q_1$) and sell ($q_2$) decisions, where $\omega_1 = \omega_2 + q_1 - q_2$. 

The set of admissible buying ($q_1$) and selling ($q_2$) decisions is defined by 
\begin{align}
\mathscr{A}(t, \omega_2) = \left\{ (q_1,q_2) \mid 0 \leq \omega_2 + q_1 - q_2 \leq \overline{U}, q_1 \geq 0, q_2 \geq 0 \right\}.
\end{align}
This means that we cannot sell more than our current inventory, and we cannot buy more than fits in our inventory capacity at the producer $\overline{U}$. Buying hydrogen comes at a unit price of $b_1$, and selling hydrogen comes at a unit price of $b_2 < b_1$. It is important to highlight that since $b_1$ and $b_2$ are fixed, speculative motives for buying or selling in the market are eliminated. Next to these costs, the producer incurs a fixed ordering cost each time the producer decides to buy from the external source ($K_1$) and a so-called fixed emergent purchase cost ($K_2$). The latter cost reflects the case when the current inventory is insufficient to meet today's customer demand, i.e. when $\omega_2 < 0$. The difference with the fixed ordering cost is that the producer may choose to keep a safety stock, to ensure that there is sufficient inventory to meet today's customer demand with a higher probability. According to these costs, the costs of taking decision $(q_1, q_2)$ is defined as: 
\begin{align}
    c(\omega_2,q_1,q_2) = K_1 \delta_1 + K_2 \delta_2 + b_1 q_1 - b_2 q_2,
\end{align}
where $\delta_1 = 1$ if $q_1 > 0$ and 0 otherwise, and $\delta_2 = 1$ if $\omega_2 < 0$ and 0 otherwise.

For states $s_1$ and $s_2$, we let $v(s_1)$ and $y(s_2)$ be the value functions that represent the expected future cost on the infinite time horizon. We assume that $\omega_1, \omega_2$, and $O_t$ are discretized 
with $\Phi_t$, the finite discrete support of $O_t$. Let $\mathbb{P}\{O_t = \phi_t\}$ denote the probability of having a net outflow quantity of $\phi_t$. Then, the transitions from state $s_1$ to $s_2$, and vice versa, is recursively defined by the following set of equations: 
\begin{align}
    y(t, \omega_2) &= \min_{(q_1,q_2) \in \mathscr{A}(s_2)} \Big\{ c(\omega_2,q_1,q_2) + \sum_{\phi_t \in \Phi_t } \mathbb{P}\left\{O_t = \phi_t \right\}   y(t + 1\mod T, \omega_2 +q_1-q_2 - \phi_t) \Big\}, \label{eq:trComb}
\end{align}
i.e.,
\begin{align}
    v(t,\omega_1) & =  \sum_{\phi_t \in \Phi_t } \mathbb{P}\left\{O_t = \phi_t \right\}   y(t, \omega_1 - \phi_t), \label{eq:tr1}\\
    y(t, \omega_2) &= \min_{(q_1,q_2) \in \mathscr{A}(s_2)} \Big\{ c(\omega_2,q_1,q_2) + v(t + 1\mod T, \omega_2+q_1-q_2) \Big\},\label{eq:tr2}
\end{align}
We employ Equations \eqref{eq:tr1} and \eqref{eq:tr2} instead of Equation \eqref{eq:trComb} since they offer greater ease with computations. In Equation \eqref{eq:tr1}, we let the value $v(t, \omega_1)$ being in state $s_1 = (t, \omega_1)$ be the expected value of transitioning to states $s_2 = (t, \omega_2)$ under the probability measure $\mathbb{P}$ based on the net outflow $O_t$. In Equation \eqref{eq:tr2}, we select buy and sell decision in order to minimize the expected future cost of the operational-level decisions. This system of equations can easily be solved by methods such as value iteration after imposing an appropriate discretization of the inventory levels and positions.

Summarizing, depending on the tactical-level inventory routing solution $\textbf{x}$ and associated with $\textbf{x}$ the transportation delivery schedule imposed by the selected clusters $\overline{\mathscr{R}}$, we define $\Pi(\mathbf{x})$ as the set of all feasible policies of the MDP. The actions associated with the minimum values of $y(t,\omega_2)$ define the optimal dynamic policy for the purchasing decisions of the producer, denoted as $\pi(\mathbf{x})$. Upon that, we define $C(\pi(\mathbf{x}))$ as the expected cycle costs of the producer's buying and selling decisions.

\subsection{The Joint Optimization Problem} \label{sect:modelJoint}
The goal of the SCIRP is to determine both a transportation delivery schedule $\mathbf{x}$ (and thus $\overline{\mathscr{R}}$) and a dynamic policy $\pi(\mathbf{x})$ for the producer's purchasing decisions that minimizes the total expected cost. The following optimization problem reflects these joint decisions:
\begin{align}
    \text{min  } & \sum_{r \in \mathscr{R}} c_r x_r  + C(\pi(\mathbf{x}))\\
    \text{s.t  } & \sum_{r \in \mathscr{R}} \beta_r^i x_r = 1 & \forall i \in \mathscr{N}, \label{custAssign}\\
    & \pi(\mathbf{x}) \in \Pi(\mathbf{x}), \label{MDP-first} \\
    & x_r \in \{0,1\} & \forall r \in \mathscr{R}.
\end{align}
Here, Constraint \ref{custAssign} ensures that the customer clusters are selected in such a way that each customer is assigned to exactly one cluster. Constraint \ref{MDP-first} defines the operational-level purchasing problem for the selected clusters. We provide an illustrative problem instance of the SCIRP and a feasible solution in Appendix \ref{sect:toy}, which details all the aforementioned elements of our problem. In the next section, we introduce a generic approach using ideas from parameterized cost function approximations that balances the transportation delivery schedules of multiple clusters and thus decrease the total expected cycle cost in the SCIRP.

\section{Parameterized Cost Function Approximation Approach} \label{sect:soln}

\begin{figure}[ht]
    \centering
    \includegraphics[width=14cm]{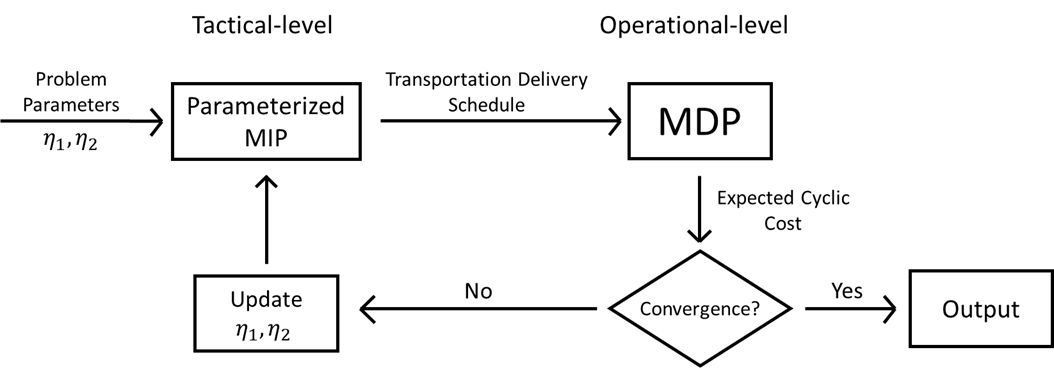}
    \caption{The framework of our approach.} \label{fig:approachFrame}
\end{figure}

The joint optimization problem is highly non-linear. The outcome of the tactical-level decision ($\textbf{x}$, the transportation delivery schedule) determines the state, action, and transition space of the MDP underlying the operational-level dynamic purchasing decisions. However, for a fixed transportation delivery schedule we can easily obtain the optimal dynamic purchasing decisions by solving the underlying MDP to optimality via value iteration. In this section, we introduce a parameterized cost function approximation (CFA) approach to handle the high complexity of the problem and improve upon solutions that ignore the relationship between the dynamic purchasing decision and the transportation delivery schedule. Our approach is based on iteratively solving a parameterized mixed-integer programming model to obtain a transportation delivery schedule and subsequently solving an MDP for the dynamic purchasing decision, see Figure \ref{fig:approachFrame}.

Overall, this approach has various benefits. Firstly, by decomposing the approach into a MIP and an MDP part, we can rely on commercial MIP solvers and value iteration, easing the computational burden of our approach. Secondly, we observe in preliminary experiments that balancing the outflow from the producer on cycle periods yields lower costs for the dynamic purchasing decision in general. For example, overlapping deliveries in the same period result in larger variances of the net outflow from the producer, resulting in more purchasing costs from the external supplier (see Appendix \ref{sect:toy}). Thus, we can steer the transportation delivery schedules by parameterizing the set-partitioning model so that the expected cost of the combined tactical and operation-level decisions is reduced. 

The remainder of this section is structured as follows.
Firstly, we will explain the idea behind our approach and present the modified 
set-partitioning model. Secondly, we discuss how we efficiently enumerate customer clusters to populate the modified set-partitioning model. Thirdly, we detail how we numerically solve the MDP that is associated with the dynamic purchasing decision. Finally, we introduce an efficient algorithm to search through the parameter space of our CFA approach by iterative solving the modified set-partitioning problem and MDP effectively.

\subsection{Modified Set-Partitioning Model}\label{sec:modmodel}

The modified set-partitioning model is subject to two parameters, $\eta_1 \geq 0$ and $\eta_2 \geq 0$, that modify the objective function by introducing a cost that depends on the combination of selected clusters and their transportation delivery schedules. The modified objective function consists of three parts. First, we minimize the expected cyclic tactical-level costs $\sum_{r \in \mathscr{R}} c_r x_r $, similar as in the set-partitioning model \eqref{initStart}-\eqref{initEnd}. Secondly, we introduce decision variables $M_{1t}$ and $M_{2t}$ and scale them with the parameters $\eta_1$ and $\eta_2$ and add this to the objective function, i.e.,  $\frac{\eta_1}{T} \sum_{t \in \mathscr{T}} M_{1t} + \frac{\eta_2}{T} \sum_{t \in \mathscr{T}}  M_{2t}$.
The decision variables penalize features of the combined selection of clusters for the transportation delivery schedule (i.e., the tactical-level decision). We select these features as the deviations between cycle periods in terms of the mean and the variance of the demand of the selected clusters (see Constraints \eqref{cfaCons1} and \eqref{cfaCons2}). The modified set-partitioning model is given by: 
\begin{align}
    \text{min  } & \sum_{r \in \mathscr{R}} c_r x_r  + \frac{\eta_1}{T} \sum_{t \in \mathscr{T}} M_{1t} + \frac{\eta_2}{T} \sum_{t \in \mathscr{T}}  M_{2t} \label{updMod2Start}\\
    \text{s.t  } & \sum_{r \in \mathscr{R}} \beta_r^i x_r = 1 & \forall i \in \mathscr{N}, \label{custAssign2}\\ 
    & - M_{1t} \leq  \frac{ \sum_{t' \in \mathscr{T}} \sum_{r \in \mathscr{R}} \Delta_r^{t'} x_r}{T} - \sum_{r \in \mathscr{R}} \Delta_r^t x_r  \leq   M_{1t} & \forall t \in \mathscr{T}, \label{cfaCons1} \\
    & - M_{2t} \leq \frac{ \sum_{t' \in \mathscr{T}} \sum_{r \in \mathscr{R}} \Lambda_r^{t'} x_r}{T} - \sum_{r \in \mathscr{R}} \Lambda_r^t x_r \leq   M_{2t} & \forall t \in \mathscr{T}, \label{cfaCons2} \\
    & x_r \in \{0,1\} & \forall r \in \mathscr{R}, \\
    & M_{1t}, M_{2t} \geq 0 & \forall t \in \mathscr{T}.  \label{updMod2End}
\end{align}

In Constraint \eqref{cfaCons1}, $\Delta_r^t$ is the expected demand at delivery period $t$ of cluster $r$, and $0$ otherwise; i.e., $\Delta_r^t \coloneqq \mathbb{E}\left[ \sum_{i \in \mathscr{N}_r}  q_r^{it} \right]$. 
Thus, $(\sum_{t \in \mathscr{T}} \sum_{r \in \mathscr{R}} \Delta_r^t x_r )/ T$ is the average demand per period over the cycle. This constraint is defined for each cycle period, and we set the deviation from the mean demand equal to $M_{1t}$ for each period. 
Similarly, Constraint \eqref{cfaCons2} aims to minimize the deviation per period of the variance of the customer demand generated by the customer clusters. Here, $\Lambda_r^t$ is the variance of the demand of cluster $r$ at delivery period $t$, i.e., 
$\Lambda_r^t \coloneqq \text{Var}\left[  \sum_{i \in \mathscr{N}_r} q_r^{it} \right]$, and $M_{2t}$ is set equal to the deviation from the average demand variance per period $t$.

In our SCIRP, the modified set-partitioning model creates a balance in the demand that the producer faces, which might lead to lower costs of buying or selling hydrogen in the dynamic purchasing decision. The decision maker can select higher values for $\eta_1$ and $\eta_2$ and the model automatically creates replenishment schemes that balance the demand over the periods in the cycle. Our experiments indeed confirm that increasing $\eta_1$ and $\eta_2$ increases the costs associated with the tactical-level while it decreases the costs of purchasing from the external supplier in general.

Finally, we like to stress our presented approach is fairly general. The characteristics we select for the SCIRP (demand and variance variability) are valid for various types of demand distributions of customers. Furthermore, initial experiments have shown that  alternative penalizations regarding, for example, the demand variances per period, lead to similar outcomes, which offers opportunities to tailor our approach more specifically in case other problem settings are studied. 
Of course, without a known demand distribution creating tactical-level decisions on its own becomes much more complex. Nevertheless, assuming a method exists to derive tactical-level decisions, our modeling strategy provides a simple but elegant way to connect this to the quality of operational-level decisions. This even surpasses the applications we study in this paper. 

\subsection{Populating the Set-Partitioning Model}

The modified set-partitioning model requires as input the customer clusters, their costs, and associated decisions such as the base-stock policy and delivery periods. In this section, we provide an algorithm that recursively finds all the feasible customer clusters that are candidates to be selected in an optimal solution. Let $\mathscr{R}$ be the set of customer groups. Each customer group $r \in \mathscr{R}$ is defined by i) a set of customers, ii) the shortest vehicle route among the customers in $r$, iii) the base-stock policy to meet the predefined service-level $\alpha$ for each customer $i \in \mathscr{N}_r$ , and iv) the delivery periods in the cycle. 

The algorithm works as follows. First, we compute all the base-stock levels required for covering demand $n$ periods ahead, $F^{-1}_{in}(\alpha)$. Next, we enumerate all potential customer clusters assuming that we deliver them at each period, i.e., a full transportation delivery schedule. Customers are recursively added to clusters until Constraint \eqref{cc2} is violated. For each feasible customer cluster assuming a full transportation delivery schedule, the algorithm then checks the feasibility of all transportation delivery schedules and calculates the relevant parameters including costs and the penalized features in the modified set-partitioning model. The pseudo-code of the algorithm is given in Algorithm \ref{alg:setR} in Appendix \ref{appendix2}.

\subsection{Obtaining Optimal Dynamic Purchasing Decisions} \label{solDynamic}

We solve the modified set-partitioning model subject to the enumerated cluster set $\mathscr{R}$ given parameters $(\eta_1,\eta_2)$. The tactical-level decision, $\mathbf{x}$, is then input for the MDP that describes the dynamic purchasing decisions. We use value iteration \citep[see, e.g.,][]{puterman2014markov} to solve the MDP part. The obtained policy is simulated to derive the relevant cost components.

\subsection{Line Search on the Parameterized CFA}
Our final approach for solving the SCIRP iterates between solving the modified set-partitioning model with given $(\eta_1, \eta_2)$ and evaluating the associated optimal dynamic purchasing decisions. In this section, we provide an efficient search algorithm in the $(\eta_1,\eta_2)$ space to find the values of $(\eta_1, \eta_2)$ that results in the joint solution of minimum cost. We make use of the observation in preliminary experiments that the total cost function is typically unimodal in $\eta_1$ and $\eta_2$.
Using this idea, we propose a line search on the parameterized CFA where we iteratively update either $\eta_1$ or $\eta_2$ in each iteration. The algorithm, which we call Line Search, is provided in Algorithm \ref{alg:line} in Appendix \ref{appendix2}. 

The algorithm starts by initializing $(\eta_1,\eta_2) = (\epsilon,\epsilon)$, with an arbitrarily small positive value to activate the CFA approach, and then iteratively updates either $\eta_1$ or $\eta_2$. We first incrementally increase the parameter $\eta_1$ until the expected cycle cost is non-increasing, or until the predetermined upper bound of $\eta_1$ is reached. The algorithm, then, fixes $\eta_1$ and incrementally increases $\eta_2$ until the same stopping criteria are met. Three different operations are executed on each selected pair of parameters, which define a feasible solution for the corresponding pair of parameters and yield an expected cycle cost of the feasible solution. First, the MIP model of parameterized CFA is solved. Second, the corresponding MDP upon the tactical-level solution is solved. Thirdly, the optimal policy of MDP is simulated to derive the costs of the MDP policy. Among the solutions that are obtained during the process, the one with the minimum expected cycle cost is selected as the best solution.

\section{Computational Experiments} \label{sect:comp}

This section evaluates our Line Search method in terms of solution quality and computation time. The experiments are performed using an Intel Core i7-10750H CPU (2.6 GHz) with 16GB of RAM. The algorithms are implemented in C++20 in combination with Gurobi 9.1.1. We first describe the base system of our test instances in Section \ref{sect:Base}. Next, our approach is compared to a step-by-step solution in 
 Section \ref{sect:00}; where the transportation delivery schedules are obtained assuming infinite supply, thus ignoring the impact and cost of the dynamic purchasing decisions similar to the extant literature. Afterward, we compare our Line Search method with a full parameter grid search on medium-sized instances in Section \ref{sect:brute-grid}.

\subsection{Base System} \label{sect:Base}
We consider a base system given in Table \ref{baseParam}. In each instance, $15$ customers are randomly allocated into a region of $[0,10] \times [0,10]$. For each customer, we assign the mean demand in the range of $[100, 400]$ kg per period, where the mean supply equals the total mean demand. The standard deviation of customer demands is randomized in the range of $[0.025,0.05]$ times the mean demand. For the producer, the standard deviation equals $0.15$ times the mean supply. 

\begin{table}[h]
\centering
\caption{Parameter set of the base system.} \label{baseParam}
\begin{tabular}{ccccccccccccc}
\toprule $W$ & $w$ & $e$ & $Q$ & $h$ & $\alpha$ & $\gamma$ & $U$ & $C$ & $K_1$ & $K_2$ & $b_1$ & $b_2$  \\ \midrule
100 & 20 &	10 &	1000 &	0.05 &	0.95 &	0.9 &	1000 &	4500 &	3000 &	15000 &	25 &	2 \\  \bottomrule
\end{tabular}
\end{table}

We discretize the continuous sets of inventory levels to the nearest integer values. As the value iteration algorithm provides $\epsilon$-optimality, we conducted preliminary experiments and set $\epsilon = 0.1$ to get the right balance between computational efficiency and the quality of the solution. That is, the obtained solutions will not significantly improve for smaller values of $\epsilon$. After we obtain the $\epsilon$-optimal solutions of the value iteration, we simulate the output policy with $30$ million periods to derive associated costs.

\subsection{Comparison of Our Approach with the Step-By-Step Solution} \label{sect:00}

The SCIRP in previous research is studied while assuming an infinite supply. To compare our Line Search method with the literature, we solve our benchmark instances also by a step-by-step solution. This ignores the relation between the tactical and operational-level decisions, and thus solely optimizes the static, inventory-routing decisions. The resulting solution is simply evaluated in the associated MDP to obtain the cost of the optimal dynamic purchasing decision.

The results are given in Table \ref{table:00}. Upon the base case, we vary the number of customers between 10, 15, and 20, the cycle length between 7 and 10, the vehicle capacities between 800 and 1000, and demand uncertainty levels between  `L' and `H'. A demand uncertainty of $L$ is the base case, while for $H$ we randomize the customer demands in the range $[0.02, 0.1]$ times the mean demand. For each scenario, $20$ instances are generated and solved, where $\#_\text{rep}$ is the number of solutions that are solved within the set memory limit. $\Delta(\%)$ is the increase in cost of the step-by-step solution compared to our approach. For both our approach and step-by-step, the average values of the objective function values and the solution times are reported. For our Line Search approach, we also report the average number of iterations found on the search algorithm ($\#_{\text{iter}}$), where the solution times of Line Search are given as the average of each iteration of the algorithm. 

\begin{table}[h]
\centering
\caption{Performance of the Line Search compared to the step-by-step solution.}
\label{table:00}
\resizebox{\textwidth}{!}{
\begin{tabular}{lllrrrrrrrrrrrrr} \toprule
& & & & & \multicolumn{2}{l}{$\Delta(\%)$} & \multicolumn{4}{c}{\text{Step-By-Step}} & \multicolumn{5}{c}{\text{Line Search}} \\
\cmidrule(r){6-7}\cmidrule(lr){8-11}\cmidrule(l){12-16}
\multicolumn{1}{c}{$N$} &
\multicolumn{1}{c}{$T$} &
\multicolumn{1}{c}{Unc.} & \multicolumn{1}{c}{$Q$} & \multicolumn{1}{c}{$\#_\text{rep}$} & \multicolumn{1}{c}{avg.} & \multicolumn{1}{c}{max} & \multicolumn{1}{c}{$\text{obj}_\text{MIP}$} & \multicolumn{1}{c}{$\text{obj}_\text{MDP}$} & \multicolumn{1}{c}{$\text{time}_\text{MIP}$} & \multicolumn{1}{c}{$\text{time}_\text{MDP}$} & \multicolumn{1}{c}{$\text{obj}_\text{MIP}$} & \multicolumn{1}{c}{$\text{obj}_\text{MDP}$} & \multicolumn{1}{c}{$\text{time}_\text{MIP}$} & \multicolumn{1}{c}{$\text{time}_\text{MDP}$} & \multicolumn{1}{c}{$\#_\text{iter}$} \\ \midrule
10 & 7  & L & 800  & 20 & 11.0 & 48.3 & 10218 & 7642  & 0.03 & 0.14 & 10338 & 5809  & 9.5  & 0.15 & 11.3 \\
   &    &   & 1000 & 20 & 15.5 & 62.6 & 8657  & 8184  & 0.05 & 0.13 & 8776  & 5870  & 7.0  & 0.15 & 9.7  \\
   &    & H & 800  & 20 & 8.0  & 27.9 & 11216 & 7907  & 0.03 & 0.12 & 11315 & 6404  & 3.1  & 0.13 & 8.7  \\
   &    &   & 1000 & 20 & 8.8  & 44.3 & 9446  & 7966  & 0.03 & 0.13 & 9582  & 6496  & 7.2  & 0.14 & 8.8  \\
   & 10 & L & 800  & 19 & 17.4 & 51.2 & 9932  & 8716  & 0.17 & 0.16 & 10054 & 5889  & 25.7 & 0.22 & 11.5 \\
   &    &   & 1000 & 19 & 32.6 & 191.5  & 8155  & 10423 & 0.35 & 0.14 & 8302  & 5805  & 84.9 & 0.22 & 15.5 \\
   &    & H & 800  & 18 & 18.5 & 72.6 & 11327 & 10059 & 0.14 & 0.15 & 11443 & 6671  & 23.6 & 0.19 & 10.7 \\
   &    &   & 1000 & 16 & 19.5 & 82.2 & 9388  & 9761  & 0.26 & 0.15 & 9547  & 6543  & 85.8 & 0.20 & 10.0 \\
15 & 7  & L & 800  & 20 & 2.6  & 10.5 & 14942 & 8619  & 0.04 & 0.20 & 15049 & 7914  & 5.2  & 0.19 & 10.1 \\
   &    &   & 1000 & 20 & 4.6  & 12.7 & 12520 & 9103  & 0.09 & 0.19 & 12655 & 8023  & 19.1 & 0.19 & 10.0 \\
   &    & H & 800  & 20 & 2.0  & 8.2  & 16574 & 9110  & 0.05 & 0.18 & 16656 & 8530  & 5.4  & 0.17 & 8.2  \\
   &    &   & 1000 & 20 & 5.2  & 19.4 & 14010 & 9712  & 0.09 & 0.18 & 14107 & 8539  & 24.7 & 0.18 & 7.9  \\
   & 10 & L & 800  & 19 & 6.7  & 18.5 & 15096 & 9692  & 0.44 & 0.24 & 15164 & 8020  & 119.3  & 0.26 & 11.5 \\
   &    &   & 1000 & 18 & 12.3 & 40.6 & 12190 & 10604 & 1.11 & 0.21 & 12305 & 8055  & 115.6  & 0.25 & 11.6 \\
   &    & H & 800  & 16 & 5.3  & 12.2 & 16656 & 10057 & 0.32 & 0.25 & 16723 & 8656  & 109.2  & 0.24 & 10.7 \\
   &    &   & 1000 & 14 & 8.0  & 39.5 & 13601 & 10238 & 0.94 & 0.24 & 13725 & 8422  & 50.9 & 0.25 & 9.1  \\
20 & 7  & L & 800  & 20 & 1.1  & 4.2  & 19354 & 10246 & 0.11 & 0.21 & 19468 & 9822  & 38.3 & 0.19 & 9.5  \\
   &    &   & 1000 & 19 & 3.0  & 10.2 & 16251 & 10802 & 0.25 & 0.18 & 16311 & 10003 & 80.8 & 0.18 & 8.1  \\
   &    & H & 800  & 20 & 0.9  & 4.3  & 21406 & 10728 & 0.10 & 0.20 & 21480 & 10361 & 28.0 & 0.18 & 8.5  \\
   &    &   & 1000 & 20 & 2.1  & 5.6  & 17861 & 11057 & 0.20 & 0.21 & 17963 & 10356 & 98.7 & 0.19 & 7.3  \\
   & 10 & L & 800  & 16 & 4.7  & 15.1 & 19302 & 11847 & 0.99 & 0.28 & 19360 & 10384 & 644.5  & 0.27 & 12.8
\\ \bottomrule                  
\end{tabular}
}
\end{table}

We observe a sharp decrease in the cost of the dynamic purchasing decision in all solutions at the expense of only a slight increase in transportation delivery schedule costs when using Line Search compared to the step-by-step solution. As expected, the Line Search solutions are strictly better than the step-by-step solution. The average differences range between 32.6\% and 0.9\%, while the maximum difference is 191.5\%. This represents a 191.5\% increase in costs when using the step-by-step solution compared to our approach. For an increasing number of customers, we see that the differences reduce when comparing the average and the maximum cost difference between Line Search and the step-by-step solution ($\Delta(\%)$). However, the Line Search appears to be very robust compared to the step-by-step solution, indicated by the significant maximum differences for the largest instances. 
Furthermore, we observe that the solution time of the MDP does not depend on the number of customers, whereas the MIP does, because we use the same number of discretized state values in its formulation. The solution time of MIP is strongly correlated to the selected values of $\eta_1$ and $\eta_2$, which is explained further in Section \ref{sect:brute-grid}.

\subsection{Comparison of Our Approach with Grid Search} \label{sect:brute-grid}
In this section, we test the solution quality of Line Search by comparing it with a grid search on the $(\eta_1, \eta_2)$ space. For the grid search, we let $\eta_1 \in \{0, 0.0001, 1, 2, 3, 4, 5, 6, 7, 8\}$ and $\eta_2 \in \{0, 0.0001, 0.5, 1.0, 1.5, 2.0, 2.5, 3.0, 3.5, 4\}$, except for the high demand uncertainty systems. For this or when $N=20$, we let $\eta_1 \in \{0, 0.0001, 1, 2, 3, 4, 5, 6, 7\}$ and $\eta_2 \in \{0, 0.0001, 0.5, 1.0, 1.5, 2.0, 2.5\}$. We add a low value of $0.0001$ to these sets in order to test how activation of our approach with an epsilon value affects the solution quality.

\begin{figure}[h]
     \centering
     \begin{subfigure}[t]{0.48\textwidth}
         \centering
          \begin{tikzpicture}[scale=0.8]
\begin{axis}[xlabel=$\eta_1$,ylabel=$\eta_2$, zlabel = $\#_{\text{OPT}}$ , view = {190}{40} ]

\addplot3[
    surf, shader=faceted interp , samples = 2
] 
coordinates {
	(0,0,1)	(0,0.0001,1)	(0,0.5,1)	(0,1,1)	(0,1.5,1)	(0,2,2)	(0,2.5,2)	(0,3,2)	(0,3.5,2)	(0,4,2)
	
(0.0001,0,1)	(0.0001,0.0001,1)	(0.0001,0.5,1)	(0.0001,1,1)	(0.0001,1.5,1)	(0.0001,2,2)	(0.0001,2.5,2)	(0.0001,3,2)	(0.0001,3.5,2)	(0.0001,4,2)

(1,0,6)	(1,0.0001,6)	(1,0.5,7)	(1,1,4)	(1,1.5,4)	(1,2,4)	(1,2.5,4)	(1,3,4)	(1,3.5,4)	(1,4,3)

(2,0,8)	(2,0.0001,8)	(2,0.5,6)	(2,1,7)	(2,1.5,6)	(2,2,5)	(2,2.5,5)	(2,3,5)	(2,3.5,5)	(2,4,4)

(3,0,9)	(3,0.0001,9)	(3,0.5,9)	(3,1,6)	(3,1.5,7)	(3,2,6)	(3,2.5,5)	(3,3,4)	(3,3.5,4)	(3,4,4)

(4,0,5)	(4,0.0001,5)	(4,0.5,5)	(4,1,5)	(4,1.5,5)	(4,2,5)	(4,2.5,5)	(4,3,5)	(4,3.5,4)	(4,4,4)

(5,0,4)	(5,0.0001,4)	(5,0.5,4)	(5,1,4)	(5,1.5,5)	(5,2,5)	(5,2.5,5)	(5,3,5)	(5,3.5,4)	(5,4,4)

(6,0,4)	(6,0.0001,4)	(6,0.5,4)	(6,1,4)	(6,1.5,4)	(6,2,5)	(6,2.5,5)	(6,3,4)	(6,3.5,4)	(6,4,4)

(7,0,4)	(7,0.0001,4)	(7,0.5,4)	(7,1,4)	(7,1.5,4)	(7,2,5)	(7,2.5,4)	(7,3,4)	(7,3.5,4)	(7,4,4)

(8,0,3)	(8,0.0001,3)	(8,0.5,4)	(8,1,4)	(8,1.5,4)	(8,2,4)	(8,2.5,4)	(8,3,4)	(8,3.5,4)	(8,4,4)

};
\end{axis}
\end{tikzpicture}
         \caption{\footnotesize $N=15$, $T=7$, Unc = $L$, $Q=1000$.}
         \label{fig:gridOpt1}
     \end{subfigure}
     \hfill
     \begin{subfigure}[t]{0.48\textwidth}
         \centering
         \begin{tikzpicture}[scale=0.8]
\begin{axis}[xlabel=$\eta_1$,ylabel=$\eta_2$, zlabel = $\#_{\text{OPT}}$ , view = {190}{40} ]

\addplot3[
    surf, shader=faceted interp , samples = 2
] 
coordinates {
	(0,0,6)	(0,0.0001,6)	(0,0.5,4)	(0,1,3)	(0,1.5,1)	(0,2,1)	(0,2.5,1)
	
(0.0001,0,6)	(0.0001,0.0001,6)	(0.0001,0.5,4)	(0.0001,1,3)	(0.0001,1.5,1)	(0.0001,2,1)	(0.0001,2.5,1)

(1,0,10)	(1,0.0001,10)	(1,0.5,6)	(1,1,2)	(1,1.5,1)	(1,2,1)	(1,2.5,1)

(2,0,8)	(2,0.0001,8)	(2,0.5,2)	(2,1,2)	(2,1.5,1)	(2,2,1)	(2,2.5,1)

(3,0,5)	(3,0.0001,5)	(3,0.5,3)	(3,1,2)	(3,1.5,1)	(3,2,1)	(3,2.5,1)

(4,0,5)	(4,0.0001,5)	(4,0.5,3)	(4,1,1)	(4,1.5,1)	(4,2,1)	(4,2.5,1)

(5,0,4)	(5,0.0001,4)	(5,0.5,2)	(5,1,1)	(5,1.5,1)	(5,2,1)	(5,2.5,1)

(6,0,3)	(6,0.0001,3)	(6,0.5,2)	(6,1,1)	(6,1.5,1)	(6,2,1)	(6,2.5,1)

(7,0,2)	(7,0.0001,2)	(7,0.5,2)	(7,1,1)	(7,1.5,1)	(7,2,1)	(7,2.5,1)
};
\end{axis}
\end{tikzpicture}
         \caption{\footnotesize $N=20$, $T=7$, Unc = $H$, $Q=800$.}
         \label{fig:gridOpt2}
     \end{subfigure} 
        \caption{Number of best solutions on each pair of $\eta$ values.}
        \label{fig:gridOpt}
\end{figure}
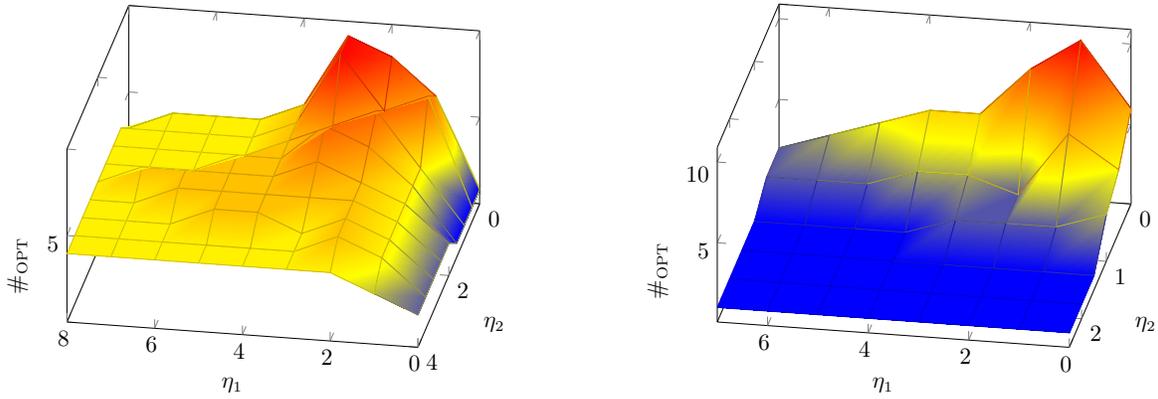

We observe that the gap between our approach and the best found in the grid search is $0.07\%$ on average for the instances in Section \ref{sect:00}, while $88\%$ of the time Line Search finds the best objective function value in the grid. In Figure \ref{fig:gridOpt}, we plot how many times the best solution is on a particular parameter pair in the grid, including all alternative solutions. 

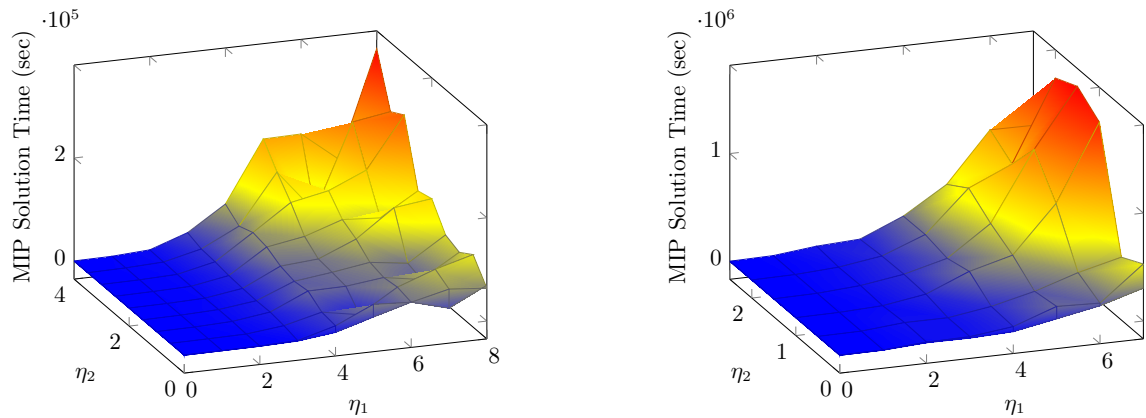
\begin{figure}[h]
     \centering
     \begin{subfigure}[t]{0.48\textwidth}
         \centering
          \begin{tikzpicture}[scale=0.8]
\begin{axis}[xlabel=$\eta_1$,ylabel=$\eta_2$, zlabel = MIP Solution Time (sec), view = {-20}{25} ]

\addplot3[
    surf, , fill= white, shader=faceted interp 
] 
coordinates {
	(0,0,78.2609666666667)	(0,0.0001,119.009366666667)	(0,0.5,362.925333333333)	(0,1,389.018666666667)	(0,1.5,436.761666666667)	(0,2,463.939)	(0,2.5,516.122666666667)	(0,3,600.634)	(0,3.5,598.32)	(0,4,765.489333333333)
	
(0.0001,0,122.776033333333)	(0.0001,0.0001,199.658333333333)	(0.0001,0.5,279.464666666667)	(0.0001,1,340.414333333333)	(0.0001,1.5,334.478666666667)	(0.0001,2,420.893)	(0.0001,2.5,513.894)	(0.0001,3,582.069666666667)	(0.0001,3.5,588.393666666667)	(0.0001,4,653.032666666667)

(1,0,487.260666666667)	(1,0.0001,679.534333333333)	(1,0.5,829.199333333333)	(1,1,883.392)	(1,1.5,1164.225)	(1,2,958.382333333333)	(1,2.5,1426.53666666667)	(1,3,1275.78933333333)	(1,3.5,1624.35966666667)	(1,4,1430.951)

(2,0,945.331)	(2,0.0001,1168.51366666667)	(2,0.5,1738.06666666667)	(2,1,3487.493)	(2,1.5,3621.37933333333)	(2,2,2931.065)	(2,2.5,3226.22333333333)	(2,3,3531.65)	(2,3.5,4697.15666666667)	(2,4,5582.77)

(3,0,2712.80566666667)	(3,0.0001,3516.418)	(3,0.5,9638.75833333333)	(3,1,13276.8166666667)	(3,1.5,12486.55)	(3,2,13629.4466666667)	(3,2.5,15674.88)	(3,3,19647.43)	(3,3.5,21743.4433333333)	(3,4,30673.6333333333)

(4,0,11265.0766666667)	(4,0.0001,23860.1366666667)	(4,0.5,27521.31)	(4,1,45747.3033333333)	(4,1.5,41840.5666666667)	(4,2,48428.7666666667)	(4,2.5,62677.8566666667)	(4,3,68066.5166666667)	(4,3.5,57730.0133333333)	(4,4,77779.4966666667)

(5,0,32704.5733333333)	(5,0.0001,57125.98)	(5,0.5,14501.8366666667)	(5,1,46661.3)	(5,1.5,50911.1833333333)	(5,2,54073.5133333333)	(5,2.5,74930.84)	(5,3,117817.3)	(5,3.5,153936.666666667)	(5,4,195406.966666667)

(6,0,53221.9533333333)	(6,0.0001,52904.3333333333)	(6,0.5,59363.5333333333)	(6,1,55342.21)	(6,1.5,61922.3566666667)	(6,2,84270.85)	(6,2.5,99218.8)	(6,3,132393.3)	(6,3.5,69703.4666666667)	(6,4,196790.4)

(7,0,33218.5)	(7,0.0001,70065.0333333333)	(7,0.5,85757.8666666667)	(7,1,43276.1)	(7,1.5,52072.3766666667)	(7,2,136591.433333333)	(7,2.5,94349.6333333333)	(7,3,163934.133333333)	(7,3.5,234442.033333333)	(7,4,139482.2)

(8,0,69412.8)	(8,0.0001,63710.3366666666)	(8,0.5,105944.8)	(8,1,87630.6)	(8,1.5,98539.5666666667)	(8,2,145945.533333333)	(8,2.5,138291.733333333)	(8,3,262504.1)	(8,3.5,246146.5)	(8,4,345409)

};
\end{axis}
\end{tikzpicture}
         \caption{\footnotesize $N=15$, $T=7$, $Q=1000$ (avg. of $3$).}
         \label{fig:gridModelT1}
     \end{subfigure}
     \hfill
     \begin{subfigure}[t]{0.48\textwidth}
         \centering
         \begin{tikzpicture}[scale=0.8]
\begin{axis}[xlabel=$\eta_1$,ylabel=$\eta_2$, zlabel = MIP Solution Time (sec), view = {-20}{25} ]

\addplot3[
    surf, , fill= white, shader=faceted interp 
] 
coordinates {
	(0,0,148.2843)	(0,0.0001,431.863333333333)	(0,0.5,1018.62133333333)	(0,1,1067.58183333333)	(0,1.5,1301.62133333333)	(0,2,1715.46666666667)	(0,2.5,1856.54466666667)
	
(0.0001,0,324.502)	(0.0001,0.0001,419.409166666667)	(0.0001,0.5,1055.75366666667)	(0.0001,1,935.967166666667)	(0.0001,1.5,1724.06233333333)	(0.0001,2,2249.309)	(0.0001,2.5,2171.98116666667)

(1,0,2337.615)	(1,0.0001,3993.875)	(1,0.5,4509.72)	(1,1,4249.38333333333)	(1,1.5,6982.56666666667)	(1,2,7212.20333333333)	(1,2.5,8829.665)

(2,0,25441.0133333333)	(2,0.0001,30595.61)	(2,0.5,38563.07)	(2,1,20845.205)	(2,1.5,33517.0383333333)	(2,2,30139.7833333333)	(2,2.5,50883.2383333333)

(3,0,23517.0433333333)	(3,0.0001,36097.8783333333)	(3,0.5,37637.185)	(3,1,82271.6833333333)	(3,1.5,92727.3883333333)	(3,2,88063.2283333333)	(3,2.5,70607.9883333333)

(4,0,37698.925)	(4,0.0001,84180.075)	(4,0.5,97303.255)	(4,1,247901.566666667)	(4,1.5,195173.371666667)	(4,2,241339.965)	(4,2.5,238824.016666667)

(5,0,111575.663333333)	(5,0.0001,225509.816666667)	(5,0.5,191679.95)	(5,1,352911.51)	(5,1.5,472417.588333333)	(5,2,626003.283333333)	(5,2.5,487132.016666667)

(6,0,187277.975)	(6,0.0001,210833.846666667)	(6,0.5,319608.52)	(6,1,785864.13)	(6,1.5,1111296.35166667)	(6,2,803460.765)	(6,2.5,947589.916666667)

(7,0,304726.54)	(7,0.0001,528241.811666667)	(7,0.5,417639.546666667)	(7,1,1494768.6)	(7,1.5,1652819.94666667)	(7,2,1556114.81666667)	(7,2.5,943898.933333333)

};
\end{axis}
\end{tikzpicture}
         \caption{\footnotesize $N=20$, $T=7$, $Q=800$  (avg. of $6$).}
         \label{fig:gridModelT2}
     \end{subfigure} 
        \caption{Average MIP solution time (in seconds) of $3$ instances on each pair of $\eta$ values.}
        \label{fig:gridModelT}
\end{figure}

The solution time of the Line Search is highly dependent on the value of $\eta_1$ and $\eta_2$. In Figure \ref{fig:gridModelT}, we provide the average times on each parameter pair for two of our instances. We observe that with an increased value of parameters, the solution time is increased sharply. This also results in a higher chance of breaking the memory limit set for the experiments.

\section{Case Study in the Northern Netherlands} \label{sect:case}

This case is based on the green hydrogen transportation in the Northern Netherlands region during the transition towards a green hydrogen based economy as a part of the project \textit{Hydrogen Energy Applications in Valley Environments for Northern Netherlands} \citep{heavenncite}. We apply our approach at two different supply locations considering expert-reviewed growth scenarios. Experts envision that growth in the region affects the number of customers, the level of supply and demand uncertainty, and the associated costs of hydrogen production and distribution \citep{newenergycoalition2020}. This case study analyses the effect of these changes on hydrogen distribution in the Northern Netherlands. 

We use the parameter sets given in Table \ref{AllParam}. These parameters are selected based on multiple resources. For the transportation costs ($W$ per replenishment and $w$ per distance), the capacities ($Q$, $U$, and $C$ per kg), and the risks of holding green hydrogen ($h$ per kg), we decide on the parameters according to the expert reviews and the stakeholders in HEAVENN. The capacities are planned to be expanded for each location over time. Although it is expected that similar vehicles will be used in the long run, the vehicle capacity will increase over years due to technological advancements resulting in more pressured tanks. For the service level probabilities ($\alpha$ and $\gamma$) and the emergency shipment cost ($e$ per kg), we refer to \cite{sonntag2021tactical} and set them as $\alpha = 0.95$, $\gamma = 0.9$, and $e = 5$. For the buy and sell prices of hydrogen per kg ($b_1$ and $b_2$ per kg), we refer to the expert discussion with partners (such as Shell, Engie, GasUnie, and Green Planet) conducted within the HEAVENN program. The current hydrogen price varies between \euro{}$10$ to \euro{}$25$ per kg. In the long run, the targeted price of hydrogen is about \euro{}$2$ per kg, which is a competitive price to its carbon-based alternatives. We incur a  difference in the buy and the sell prices which motivates the utilization of the producer's own supply. For the fixed ordering cost $K_1$ and emergent purchase $K_2$, we consider the future projections of the region and assume that the more hydrogen is available, the less fixed cost is imposed on each purchase. 

\begin{table}[h]
\centering
\caption{Parameters for the case study.}
\label{AllParam}
\begin{tabular}{l|rrrrrrrrrrrrr} \toprule
Case & $W$   & $w$  & $e$  & $Q$    & $h$   & $\alpha$ & $\gamma$ & $U$ & $C$ & $K_1$  & $K_2$  & $b_1$ & $b_2$  \\ \midrule
Emmen & 100 & 1 & 5 & 400 & 0.05 & 0.95  & 0.9   & 200 & 500 & 1000 & 1000 & 10.0  & 2     \\
Eemshaven 1 & 100 & 1 & 5 & 700 & 0.05 & 0.95  & 0.9   & 800 & 4000 & 750 & 750 & 7.5  & 2    \\
Eemshaven 2 & 100 & 1 & 5 & 1100 & 0.04 & 0.95  & 0.9   & 1500 & 8000 & 200 & 200 & 3.5  & 2 \\ \bottomrule
\end{tabular}
\end{table}

The stakeholders in HEAVENN foresee two major producers in the following years in the Northern Netherlands region; one in Emmen, and one in Eemshaven. Emmen is planned as a producer location supplying hydrogen to the nearby industry park and the refueling stations. A global supply chain of hydrogen is envisioned in the longer term where imports and exports of hydrogen become prominent bringing forward the role of seaports. In that respect, Eemshaven is being considered for taking the lead on supplying hydrogen \citep{newenergycoalition2020}. Thus, we analyze these two regions in our case study. The customer locations are selected from the existing and potential hydrogen refueling stations (HRS) in the region. The list of locations, and corresponding mean and standard deviation values are given in Table \ref{AllDist} in Appendix \ref{appendix1}. The mean and variances of supply and demand are derived from the short- and mid-term projections within HEAVENN \citep{hiddeMSC, feikoMSC, frankMSC}.  

We solve three base cases in total in the following, one for supply in Emmen (in Section \ref{sec:Emmen}) and two for supply in Eemshaven (in Section \ref{sec:Eemshaven}). For these base cases, the expected supply per day is set equal to the expected total demand per day. We further relax this assumption in Section \ref{sec:Impacts}, where the impacts of uncertainty levels are also considered. 

\subsection{Supply from Emmen} \label{sec:Emmen}

\begin{figure}[h]
    \centering
    \includegraphics[scale=0.4]{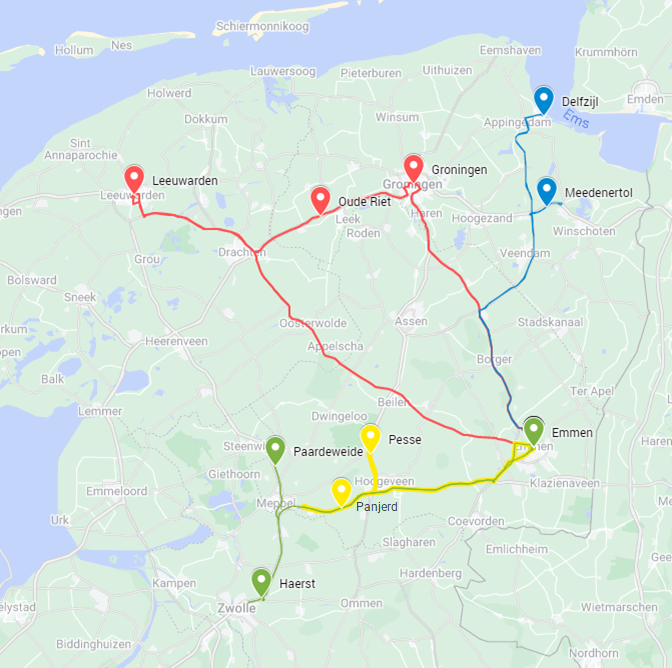}
    \caption{The routing of case of Emmen in the suggested policy.}
    \label{fig:Opt2030}
\end{figure}

Our results show a total weekly estimated cost of \euro{}$8582$ for the Emmen case. The total cost for delivering and inventory holding is found to be \euro{}$2.45$ per kg. This includes a transportation cost of \euro{}$1.15$ per kg. This seems too high once we consider the target long run selling price of \euro{}$2$ per kg, but it is reasonable with the current selling prices of \euro{}$10$ to \euro{}$25$ per kg. We expect the cost to decrease sharply within the maturation of the economy making the target price attainable, which is also validated by our case studies. The costs associated due to interaction with the external supplier (i.e. purchasing costs) is \euro{}$1.22$ per shipped kg. We observe a high reliance on the external hydrogen market due to the high supply and demand uncertainty under low volumes of the early transition states. The vehicle routes are shown in Figure \ref{fig:Opt2030}.

These routes are implemented with the transportation delivery schedule given in Table \ref{2030Schedules}. In this table, the normal distribution data of total quantity sent is given as $(\mu,\sigma^2)$ for the delivery days, where blank represents no planned delivery on that day. We notice a relatively big cluster $1$ (including major demand locations such as Groningen), where the vehicle capacity is only enough to deliver these customers on each day of the week. If a day is skipped, the targeted service level of vehicle utilization would not be met. Moreover, Groningen and Leeuwarden have relatively restricted inventory capacities (i.e., a high ratio of $\mu_i/U$), which leads to delivering these customers as frequently as possible, in order to have $\alpha = 95\%$ chance of service level for customer inventory. For the remaining clusters, the model groups two or three customers into one cluster and delivers them as seldom as possible to achieve the targeted service level. This is cost-efficient due to the relatively high transportation cost compared to the inventory holding costs. Finally, we observe that the emergency shipment is almost never used, only $0.003\%$ of the total hydrogen is expected to be delivered by additional vehicles. Overall, we see the transportation and the purchasing costs incur $96\%$ of the total costs in this case with the remaining amount barring
the costs for inventory holding and emergency shipment.

\begin{table}[h]
    \centering
    \caption{The distribution data $(\mu,\sigma^2)$ of weekly transportation delivery schedule of Emmen case.}
    \label{2030Schedules}
    \begin{tabular}{l|ccccccc} \toprule
    Cluster & Mo & Tu & We & Th & Fr & Sa & Su \\ \midrule
 1 & $(295,35^2)$ & $(295,35^2)$ & $(295,35^2)$ & $(295,35^2)$ & $(295,35^2)$ & $(295,35^2)$ &  $(295,35^2)$ \\ 
 2 &  &  & $(184,37^2)$ &  &  & $(271,32^2)$ &   \\ 
 3 &  & $(160,27^2)$ &  & $(254,27^2)$ &   &  & $(146,33^2)$  \\ 
 4 & $(251,29^2)$ &  &  &  & $(169,34^2)$ &  &    \\ \bottomrule
    \end{tabular}
\end{table}

\subsection{Supply from Eemshaven} \label{sec:Eemshaven}

\begin{figure}[h]
     \centering
     \captionsetup[subfigure]{justification=centering}
     \begin{subfigure}[b]{0.49\textwidth}
         \centering
         \includegraphics[scale=0.4]{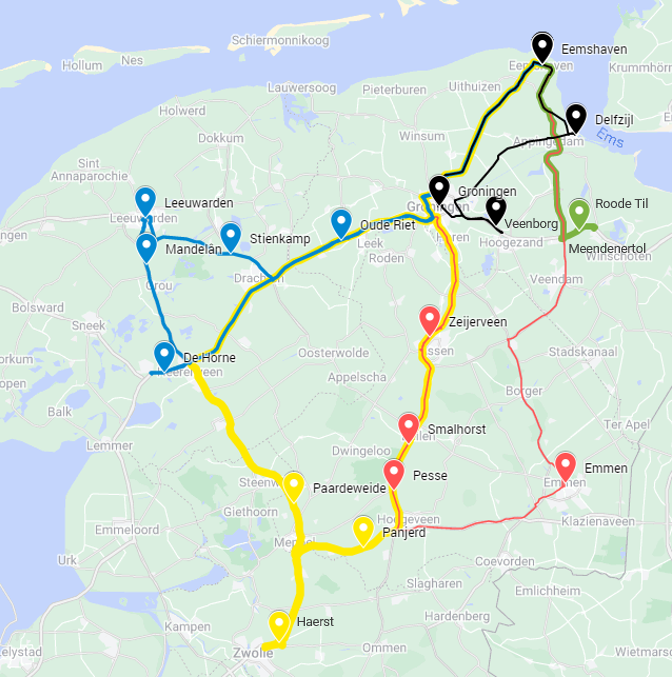}
         \caption{\footnotesize Eemshaven $1$.}
         \label{fig:Opt2040}
     \end{subfigure}
     \hfill
     \begin{subfigure}[b]{0.49\textwidth}
         \centering
         \includegraphics[scale=0.4]{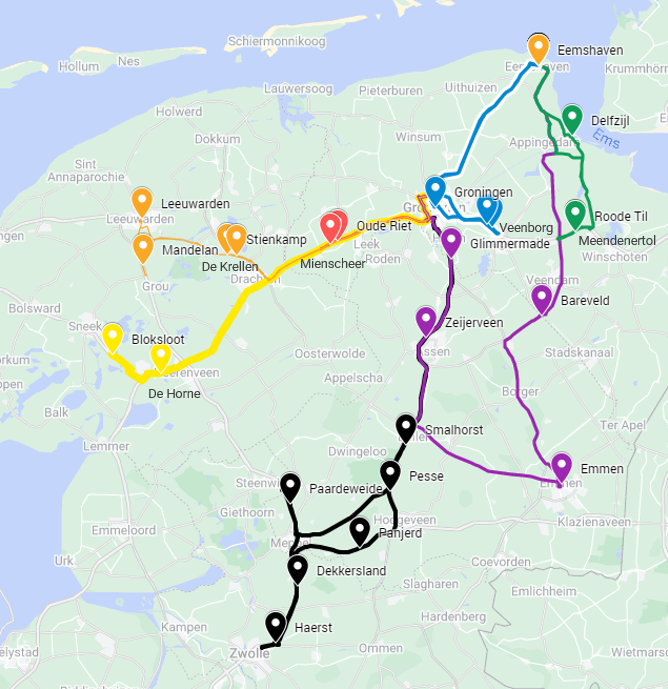}
         \caption{\footnotesize Eemshaven $2$.}
          \label{fig:Opt2050}
     \end{subfigure} 
        \caption{The routing on Eemshaven cases in the suggested policy.}
        \label{fig:Opt20402050}
\end{figure}

In the mid- and long-term for the region, Eemshaven is foreseen to be the major supplier in the Northern Netherlands. Additionally, a growth is expected in terms of number of locations and volumes of hydrogen per location compared to the short-term future, while the uncertainties decrease with a maturing hydrogen economy. We take these into account and create the instance set shown in Table \ref{AllDist} in Appendix \ref{appendix1}. The suggested weekly transportation delivery schedules are provided in Table \ref{2040Schedules} and the delivery routing is given in Figure \ref{fig:Opt2040}.

\begin{table}[h]
    \centering
    \caption{The distribution data $(\mu,\sigma^2)$ of weekly transportation delivery schedule of Eemshaven 1 case.}
    \label{2040Schedules}
    \begin{tabular}{l|ccccccc} \toprule
    Cluster & Mo & Tu & We & Th & Fr & Sa & Su \\ \midrule
 1 &  & $(600,73^2)$ &  & $(600,73^2)$ &  & $(239,73^2)$ &  $(661,52^2)$ \\ 
 2 & $(323,65^2)$ &  & $(350,53^2)$ &  & $(552,53^2)$ & &  \\
 3 & $(570,72^2)$ & $(570,72^2)$ & $(570,72^2)$ & $(570,72^2)$ & $(570,72^2)$ & $(570,72^2)$ &  $(570,72^2)$ \\ 
 4 & $(600,72^2)$ & $(600,72^2)$ & $(600,72^2)$ & $(600,72^2)$ & $(600,72^2)$ & $(600,72^2)$ &  $(600,72^2)$ \\ 
 5 & $(550,65^2)$ & $(550,65^2)$ & $(550,65^2)$ & $(550,65^2)$ & $(550,65^2)$ & $(550,65^2)$ &  $(550,65^2)$ \\  \bottomrule
    \end{tabular}
\end{table}

Since transportation costs are the major cost element of the system, we observe that the number of customers in a group increases within a higher density of customer locations in the area. Once the farthest customer from the supplier is delivered, it is more cost-efficient to deliver to other customers that are close to the farthest customer. This leads to more frequent deliveries in order to attain the vehicle capacity service level for the groups with many customers. Finally, this also minimizes the effects of customer uncertainties since combined groups have less uncertainty in total. In this case, the total weekly expected cost is \euro{}$11814$. The transportation cost per kg is \euro{}$0.54$. Note that this is less than half of the Emmen case due to the increased vehicle capacity and higher number of customers in the region which results in less traveling in between consecutive deliveries. The purchasing cost per shipped kg is \euro{}$0.16$, which is only $14\%$ of the Emmen case. This is because, with frequent deliveries and more customers on a route, the dynamic purchasing problem can provide policies with much lower costs due to reduced variances between week days. The percentage of emergency shipments to the customers is $0.27\%$ and therefore the system is mostly run on its own vehicle fleet.

\textit{\textbf{Eemshaven $2$:}} 
We study an additional Eemshaven scenario with mature hydrogen economy settings since the region is expected to become a major hub port in the long run. In this scenario, the producer is expected to replenish a total of $24$ customers. The mean levels for supply and demand are further increased. Similar to the first Eemshaven case (i.e. Eemshaven $1$), we assume a reduced uncertainty on both supply and demand in the mature stage (i.e. Eemshaven $2$), for example, the variance of daily demand of candidate HRS does not change from Eemshaven $1$ to Eemshaven $2$ while the mean value increases (see Table \ref{AllDist} in Appendix \ref{appendix1}). The solution has an expected cost of \euro{}$12877$ per week. The routes are shown in Figure \ref{fig:Opt2050}, and the transportation delivery schedule is provided in Table \ref{2050Schedules}.

\begin{table}[h]
    \centering
    \caption{The distribution data $(\mu,\sigma^2)$ of weekly transportation delivery schedule of Eemshaven 2 case.}
    \label{2050Schedules}
    \begin{tabular}{l|ccccccc} \toprule
    Cluster & Mo & Tu & We & Th & Fr & Sa & Su \\ \midrule
 1 & $(600,60^2)$ &  & $(931,60^2)$ &  &  & $(569,73^2)$ &  \\ 
 2 &  & $(931,60^2)$ &  &  & $(569,73^2)$ & & $(600,60^2)$ \\ 
 3 & $(1000,77^2)$ & $(1000,77^2)$ & $(1000,77^2)$ & $(1000,77^2)$ & $(1000,77^2)$ &  $(1000,77^2)$ &  $(1000,77^2)$ \\
 4 & $(900,85^2)$ & $(900,85^2)$ & $(900,85^2)$ & $(900,85^2)$ & $(900,85^2)$ & $(900,85^2)$ &  $(900,85^2)$ \\ 
 5 & $(850,79^2)$ & $(850,79^2)$ & $(850,79^2)$ & $(850,79^2)$ & $(850,79^2)$ & $(850,79^2)$ &  $(850,79^2)$ \\ 
 6 & $(840,69^2)$ &  & $(363,69^2)$ & $(897,49^2)$ &  & $(840,69^2)$ &   \\ 
 7 & $(900,73^2)$ & $(900,73^2)$ & $(900,73^2)$ & $(900,73^2)$ & $(900,73^2)$ & $(900,73^2)$ &  $(900,73^2)$ \\   \bottomrule
    \end{tabular}
\end{table}

In the results, we observe a similar structure to the Eemshaven $1$ case, where more and more local customers are grouped in a cluster. It is cost-efficient to have larger groups that are delivered once every day. The transportation cost per kg is \euro{}$0.33$. Even though this is the lowest compared to all cases, it is still the main cost element of the system. The purchasing cost per shipped kg is decreased to \euro{}$0.02$. In all problem instances, with less uncertainty and more capacities, the daily purchasing problem yields relatively less costs. Therefore, in a mature hydrogen market where the uncertainties are quite low, we suggest that decision makers first optimize the tactical-level operations regarding the supply of hydrogen, and then make the necessary postponements in individual transportation delivery schedules of clusters to further optimize their daily purchasing problem of distributing the hydrogen to end users as a second goal.

A summary table of each case is given in Table \ref{summaryTable} for comparison, where the costs are in \euro{}. $c^{ \textsc{t}}$ represents the transportation cost per unit, $c^{ \textsc{p}}$ is the purchasing costs per unit, $\%^\textsc{e}$ is the expected percentage of units that are transported via an emergency shipment, and lastly the total expected cost of a cycle is given.
\begin{table}[h]
\centering
\caption{The solution information of each cases.} \label{summaryTable}
\begin{tabular}{l|rrrr} \toprule
Case        & $c^{ \textsc{t}}$ & $c^{ \textsc{p}}$ & $\%^\textsc{e}$ & Total \\ \midrule
Emmen       & $1.15$            & $1.22$                    & $0.01\%$         & $8582$            \\
Eemshaven 1 & $0.54$            & $0.16$                    & $0.27\%$          & $11814$           \\
Eemshaven 2 & $0.33$            & $0.02$                    & $0.08\%$          & $12877$   \\  \bottomrule      
\end{tabular}
\end{table}

\subsection{Impacts} \label{sec:Impacts}

\begin{figure}[h]
     \centering
     \begin{subfigure}[b]{0.48\textwidth}
         \centering
         \begin{tikzpicture}[scale=0.8]
            \begin{axis}[
            xlabel style={align=center}, xlabel=Supply mean multiplier ($m_s$) \\ supply $<$ demand \hspace{1cm} supply $>$ demand,
            ylabel=Costs,
            legend pos=north east]
            \addplot[mark=*,blue] plot coordinates {
                (0.9,7529.13)(0.91,7232.14)(0.92,6945.85)(0.93,6662.98)(0.94,6397.19)(0.95,6132.47)(0.96,5868.12)(0.97,5613.23)(0.98,5369.94)(0.99,5120.79)(1,4893.81)(1.01,4663.45)(1.02,4432.52)(1.03,4215.99)(1.04,4003.61)(1.05,3792.92)(1.06,3592.97)(1.07,3396.85)(1.08,3200.73)(1.09,3017.23)(1.1,2836.53)
};
            \addlegendentry{Emmen}
    
            \addplot[color=red,mark=x] plot coordinates {
            	(0.9,12534.7)(0.91,11410)(0.92,10302.9)(0.93,9208.57)(0.94,8159.19)(0.95,7124.86)(0.96,6111.23)(0.97,5139.65)(0.98,4225.13)(0.99,3347.56)(1,2553.18)(1.01,1838.67)(1.02,1138.23)(1.03,503.95)(1.04,-75.29)(1.05,-617.78)(1.06,-1104.31)(1.07,-1564.48)(1.08,-2000.27)(1.09,-2401.44)(1.1,-2781.12)
            };
            \addlegendentry{Eemshaven 1}
            
            \addplot[color=green,mark=square] plot coordinates {
                (0.9,11612.8)(0.91,10481.2)(0.92,9328.73)(0.93,8176.09)(0.94,7067.83)(0.95,5936.38)(0.96,4810.44)(0.97,3729.9)(0.98,2666.17)(0.99,1631.35)(1,678.88)(1.01,-197.71)(1.02,-1031.04)(1.03,-1801.07)(1.04,-2511.66)(1.05,-3213.76)(1.06,-3893.15)(1.07,-4547.74)(1.08,-5216.9)(1.09,-5877.51)(1.1,-6534.12)
            };
            \addlegendentry{Eemshaven 2}
            
            \end{axis}
         \end{tikzpicture}
         \caption{\footnotesize Purchasing costs of each case.}
         \label{impact1-2}
     \end{subfigure}
     \hfill
     \begin{subfigure}[b]{0.48\textwidth}
         \centering
         \begin{tikzpicture}[scale=0.8]
            \begin{axis}[
            xlabel style={align=center}, xlabel=Supply uncertainty multiplier ($m_p$)\\ less variance \hspace{1cm} more variance,
            ylabel=Cost ratio,
            legend pos=north west]
            \addplot[mark=*,blue] plot coordinates {
                (0,0.0761)(0.1,0.0943)(0.2,0.1446)(0.3,0.2183)(0.4,0.3062)(0.5,0.4073)(0.6,0.5164)(0.7,0.6365)(0.8,0.7551)(0.9,0.8775)(1,1)(1.1,1.1174)(1.2,1.2496)(1.3,1.3711)(1.4,1.4912)(1.5,1.6086)(1.6,1.721)(1.7,1.8293)(1.8,1.934)(1.9,2.0363)(2,2.1318)

};
            \addlegendentry{Emmen}
    
            \addplot[color=red,mark=x] plot coordinates {
                (0,0.0572)(0.1,0.0671)(0.2,0.0975)(0.3,0.1475)(0.4,0.22)(0.5,0.3093)(0.6,0.4161)(0.7,0.5438)(0.8,0.68)(0.9,0.8316)(1,1)(1.1,1.1735)(1.2,1.3592)(1.3,1.561)(1.4,1.7696)(1.5,1.9792)(1.6,2.1968)(1.7,2.4174)(1.8,2.6467)(1.9,2.879)(2,3.1131)

            };
            \addlegendentry{Eemshaven 1}
            
            \addplot[color=green,mark=square] plot coordinates {
                (0,0.0386)(0.1,0.0505)(0.2,0.0793)(0.3,0.1321)(0.4,0.2053)(0.5,0.2906)(0.6,0.4081)(0.7,0.5276)(0.8,0.6687)(0.9,0.8216)(1,1)(1.1,1.1877)(1.2,1.3906)(1.3,1.5952)(1.4,1.8092)(1.5,2.0546)(1.6,2.2828)(1.7,2.5406)(1.8,2.7896)(1.9,3.0472)(2,3.3234)

            };
            \addlegendentry{Eemshaven 2}
            
            \end{axis}
         \end{tikzpicture}
         \caption{\footnotesize Percentage of purchasing costs of each case.}
         \label{fig:impact3}
     \end{subfigure} 
        \caption{Effects of the supply distribution.}
        \label{fig:impact1-3}
\end{figure}
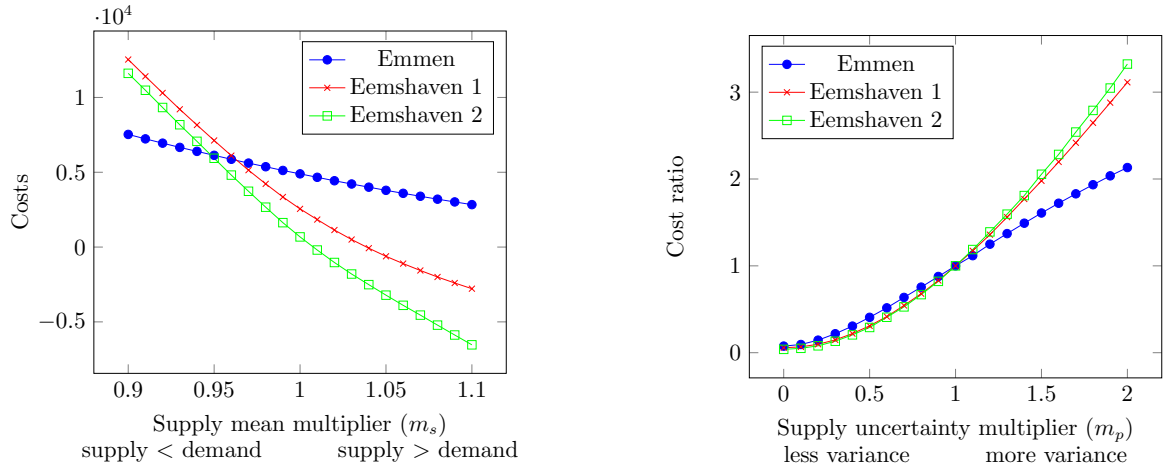

\textit{Impact of supply quantity:} For the base cases of Emmen, Eemshaven $1$, and $2$, we assume that the mean daily supply equals the mean total daily demand. However, for example, due to uncertainty in how the green hydrogen economy will emerge in society or the disruptions in supply chains, we may observe cases where this equality does not hold. In order to find the effects of the possible imbalances, we analyze all three cases by increasing and decreasing the supply amounts. Let $\mu_p$ and $\sigma_p$ be the mean and the standard deviation of daily supply, respectively. For $m_s$ in $\{0.90, 0.91, 0.92, \dots, 1.09, 1.10\}$, we analyze varying levels of supply by setting the mean equal to $m_s\mu_p$ and standard deviation equal to $m_s\sigma_p$. In Figure \ref{impact1-2}, we give the percentage change for the purchasing costs for each multiplier.

In all three cases, we observe a decrease in costs with an increase in supply. This is because increasing supply decreases the need for buying hydrogen from the external supplier, and moreover increases profits by selling the hydrogen to the external supplier. The effects of the multiplier are higher in Eemshaven cases compared to the Emmen case, possibly due to the higher supply/demand quantities which result in more purchases within, say, $1\%$ change in $m_s$. While this is an advantage for the Eemshaven cases when $m_s > 1$, it results in a sharp increase in costs if the supply is disrupted ($m_s < 1$). Thus, the more the ratio of supply/demand is decreased, the more we recommend to the decision makers to focus on their daily purchasing problem. This recommendation gains importance within the maturation of the hydrogen economy due to the higher supply and demand quantities.

\textit{Impact of supply uncertainty:} In our cases, we assume some level of uncertainty in both supply and demand. While this is decided within the expert discussions, we may observe a higher or a lower level of uncertainty in the future. We test various supply uncertainty levels to see how the system performs under different uncertainties. For $m_p \in \{0, 0.1,\dots, 1.9, 2\}$, we analyze varying levels of supply uncertainty by setting the standard deviation equal to $m_p\sigma_p$. The weekly expected purchasing costs are provided in Figure \ref{fig:impact3}. 

We observe almost a parabolic increase in purchasing costs within the increase in supply variance in Figure \ref{fig:impact3}. The purchasing costs are converging almost to zero in all cases when $m_p$ is converging to zero. Even though the system still has uncertainty in demand with $m_p = 0$, the daily purchasing problem handles that uncertainty with almost no purchase from the external supplier. This shows that, as one of our key contributions in this paper, supply uncertainty is a major factor affecting the dynamic purchasing costs, and is crucial to consider.

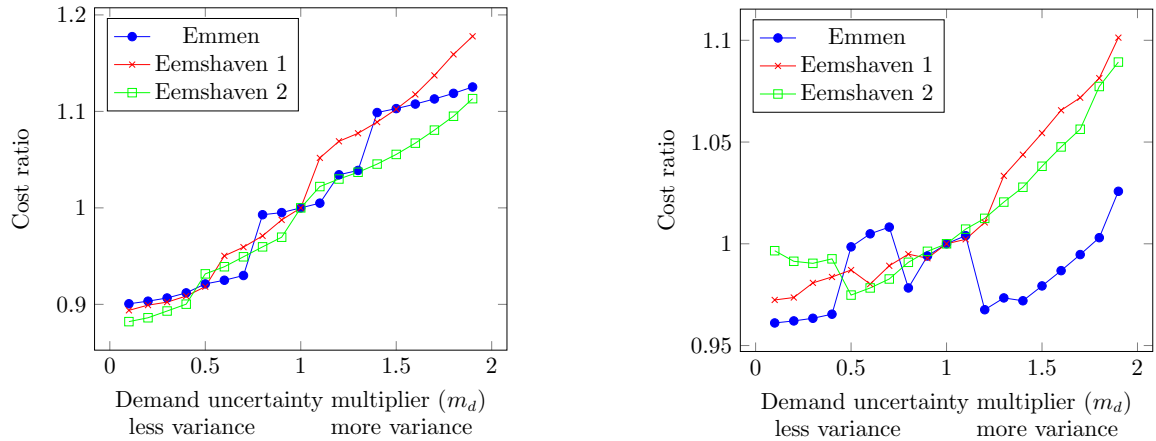
\begin{figure}[h]
     \centering
     \begin{subfigure}[b]{0.48\textwidth}
         \centering
         \begin{tikzpicture}[scale=0.8]
            \begin{axis}[
            xlabel style={align=center}, xlabel=Demand uncertainty multiplier ($m_d$) \\ less variance \hspace{1cm} more variance,
            ylabel=Cost ratio,
            legend pos=north west]
            \addplot[mark=*,blue] plot coordinates {
                (0.1,0.9006)(0.2,0.9033)(0.3,0.9066)(0.4,0.9119)(0.5,0.9212)(0.6,0.925)(0.7,0.9299)(0.8,0.9929)(0.9,0.9951)(1,1)(1.1,1.0049)(1.2,1.0343)(1.3,1.0388)(1.4,1.0988)(1.5,1.1029)(1.6,1.1076)(1.7,1.1129)(1.8,1.1187)(1.9,1.1253)

};
            \addlegendentry{Emmen}
    
            \addplot[color=red,mark=x] plot coordinates {
                (0.1,0.8938)(0.2,0.8992)(0.3,0.9023)(0.4,0.9086)(0.5,0.9183)(0.6,0.9503)(0.7,0.9593)(0.8,0.971)(0.9,0.9875)(1,1)(1.1,1.0518)(1.2,1.069)(1.3,1.0773)(1.4,1.0889)(1.5,1.1023)(1.6,1.1176)(1.7,1.1373)(1.8,1.1591)(1.9,1.1778)

            };
            \addlegendentry{Eemshaven 1}
            
            \addplot[color=green,mark=square] plot coordinates {
                (0.1,0.8821)(0.2,0.8861)(0.3,0.8931)(0.4,0.9002)(0.5,0.9316)(0.6,0.939)(0.7,0.9493)(0.8,0.9596)(0.9,0.9696)(1,1)(1.1,1.0221)(1.2,1.0299)(1.3,1.0369)(1.4,1.0454)(1.5,1.0554)(1.6,1.0671)(1.7,1.0806)(1.8,1.095)(1.9,1.1132)

            };
            \addlegendentry{Eemshaven 2}
            
            \end{axis}
         \end{tikzpicture}
         \caption{\footnotesize Tactical-level costs (transportation + inventory holding + emergency shipment).}
         \label{fig:impact2-first}
     \end{subfigure}
     \hfill
     \begin{subfigure}[b]{0.48\textwidth}
         \centering
         \begin{tikzpicture}[scale=0.8]
            \begin{axis}[
            xlabel style={align=center}, xlabel=Demand uncertainty multiplier ($m_d$) \\ less variance \hspace{1cm} more variance,
            ylabel=Cost ratio,
            legend pos=north west
            ]
            \addplot[mark=*,blue] plot coordinates {
                (0.1,0.9611)(0.2,0.9621)(0.3,0.9634)(0.4,0.9654)(0.5,0.9985)(0.6,1.0049)(0.7,1.0082)(0.8,0.9783)(0.9,0.994)(1,1)(1.1,1.004)(1.2,0.9676)(1.3,0.9734)(1.4,0.972)(1.5,0.9793)(1.6,0.9868)(1.7,0.9947)(1.8,1.003)(1.9,1.0258)

            };
           \addlegendentry{Emmen}
    
            \addplot[color=red,mark=x] plot coordinates {
                (0.1,0.9724)(0.2,0.9736)(0.3,0.9808)(0.4,0.9837)(0.5,0.9871)(0.6,0.9802)(0.7,0.9892)(0.8,0.9949)(0.9,0.9929)(1,1)(1.1,1.0022)(1.2,1.0105)(1.3,1.0334)(1.4,1.0438)(1.5,1.0544)(1.6,1.0656)(1.7,1.0718)(1.8,1.0814)(1.9,1.1013)

            };
           \addlegendentry{Eemshaven 1}
            
            \addplot[color=green,mark=square] coordinates {
                (0.1,0.9966)(0.2,0.9914)(0.3,0.9904)(0.4,0.9926)(0.5,0.9748)(0.6,0.9783)(0.7,0.9827)(0.8,0.9909)(0.9,0.9963)(1,1)(1.1,1.0072)(1.2,1.0125)(1.3,1.0205)(1.4,1.0278)(1.5,1.0381)(1.6,1.0476)(1.7,1.0563)(1.8,1.0773)(1.9,1.0893)

            };
            \addlegendentry{Eemshaven 2}
            
            \end{axis}
         \end{tikzpicture}
         \caption{\footnotesize Purchasing costs of each case. \\ }
         \label{fig:impact2-mdp}
     \end{subfigure} 
        \caption{Comparison of tactical-level and purchasing costs to the base case over demand uncertainty multiplier ($m_d$).}
        \label{fig:impact2}
\end{figure}

\textit{Impact of demand uncertainty:} As we did with supply uncertainty, we further study the effects of demand uncertainty. For $m_d$ in $\{0.1, 0.2,\dots, 1.8, 1.9\}$, we analyze varying levels of demand uncertainty by setting the standard deviation of the customer demand $i$, $\sigma_i$, equal to $m_d\sigma_i$ for all customers. The tactical-level costs and the purchasing costs are shown in Figures \ref{fig:impact2-first} and \ref{fig:impact2-mdp} for each multiplier, respectively. In both figures, the costs are given relative to the base cases, where a cost ratio of $1$ represents the same cost as the base case. 

We observe that in Figures \ref{fig:impact2-first} and \ref{fig:impact2-mdp}, the change in costs is not following a smooth function unlike the previous analysis. In the subparts where the costs are smoothly increasing, for example from $m_d = 1.4$ to $m_d = 1.8$ in Emmen case, we observe no change in the solution structure; the routing and the transportation delivery schedule. The costs are gradually increasing due to higher safety stocks and higher chances of emergency shipment. The same does not apply, however, from $m_d = 0.7$ to $m_d = 0.8$ in Emmen case, where a different solution is selected with a different route and associated daily purchasing problem. This is either because the solution of $m_d = 0.7$ is less efficient with increased uncertainty levels, or because this solution is no longer feasible due to service target levels. Additionally, the more customers we have in a region, the more combinations of feasible clusters exist. Thus, the bumpiness in costs decreases towards a mature market, due to higher possibilities for distribution. Overall, especially in the early transition stages, we recommend decision makers to update their system solutions even with a tiny change in uncertainties in customer demand.

\section{Conclusions} \label{sect:conc}
In this study, we present a novel solution framework for jointly addressing tactical-level and operational-level decision making in the context of green-hydrogen logistics. We consider a stochastic cyclic inventory routing problem that generalizes several existing problems in the context of inventory routing. The particular innovation is the consideration of supply uncertainty, due to uncertain hydrogen production from renewable sources, in combination with cyclic transportation delivery schedules that impose deliveries at fixed periods. 

Our approach combines two models: a Mixed Integer Programming (MIP) model that solves tactical-level decisions to generate a transportation delivery schedule for hydrogen, and a Markov decision process (MDP) model that optimizes operational-level decisions for buying and selling hydrogen on the market to ensure the feasibility of the transportation delivery schedule. We propose a parameterized cost function approximation approach for the tactical-level MIP model to anticipate the cost of operational-level dynamic decisions. Our approach includes an efficient search algorithm that iteratively solves MIPs and MDPs to find the optimal parameters for the modified MIP model. This allows us to effectively balance tactical and operational considerations to maximize the overall efficiency of green-hydrogen logistics operations.

We show  that our approach provides high-quality solutions on a stylized benchmark set. First, our approach outperforms an algorithm that takes a classic two-stage approach in which first the transportation delivery schedules are created (i.e., the tactical-level decision), and afterward, the dynamic purchasing decisions are made (i.e., the operational-level decisions). In addition, our efficient search over the parameter space, as required for the cost function approximation, shows excellent performance in comparison to a grid search over the parameter space.

Furthermore, we show the results of a case study to give insights into efficient green hydrogen distribution as a part of the project \textit{Hydrogen Energy Applications in Valley Environments for Northern Netherlands} \citep{heavenncite}. We solve three expert-reviewed base scenarios. Furthermore, we analyze the impact of  various problem elements upon these base scenarios and it appears that supply uncertainty is a key factor for operational-level decisions and should be anticipated while designing transportation delivery schedules on a tactical-level.

The opportunities for future research are numerous. Some of these require no fundamental change in the model and approach, only requiring generating customer clusters in a different way. Examples include the consideration of a heterogeneous vehicle fleet, enabling the use of multiple vehicles per replenishment, and demand and supply distributions without autocorrelation but that differ per period. 
Furthermore, the problem may be dealt with in a dynamic setting, enabling nonstationary and/or auto-correlated customer demands, and the daily reoptimization of vehicle routes. This latter element would also require predicting future hydrogen production based on weather forecasts, which we deem a very interesting avenue for further research.

\ACKNOWLEDGMENT{This project has received funding from the Fuel Cells and Hydrogen 2 Joint Undertaking (now Clean Hydrogen Partnership) under Grant Agreement No 875090. This Joint Undertaking receives support from the European Union's Horizon 2020 research and innovation programme, Hydrogen Europe and Hydrogen Europe Research. Albert H. Schrotenboer has received support from the Dutch Science Foundation (NWO) through grant VI.Veni.211E.043. We thank the Center for Information Technology of the University of Groningen for their support and for providing access to the Peregrine high performance computing cluster.}












\bibliographystyle{informs2014trsc}
\bibliography{main} 

\newpage
\begin{APPENDICES} 
\section{An Example Solution} \label{sect:toy}
We discuss an example instance in this section to gain insight into the dependency between the tactical-level and operational-level decisions. The location and corresponding distances of $3$ customers are given in Figure \ref{illusLocations}, where $0$ represents the producer location. The daily demands of customers $1-3$ are normally distributed with mean values of $100,250,500$ and standard deviations of $25,50,75$, respectively. The daily supply of the producer is also normally distributed with a mean of $850$ and a standard deviation of $120$.

\begin{figure}[h]
     \centering
     \begin{subfigure}[b]{0.45\textwidth}
        \centering
        \includegraphics[width=5cm]{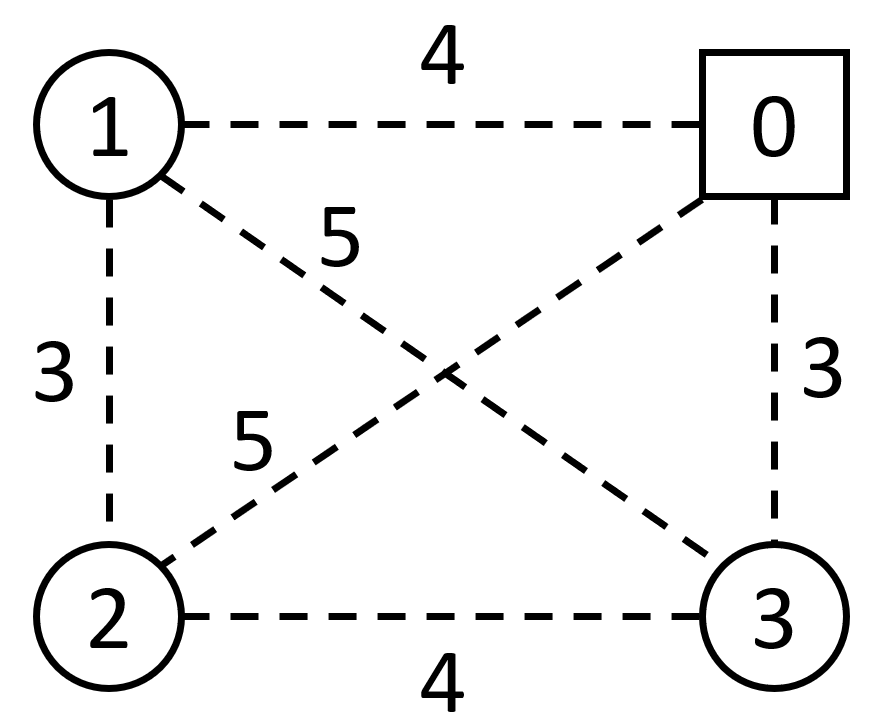}
        \caption{\footnotesize The locations and distances of example instance.}
        \label{illusLocations}
     \end{subfigure}
     \hfill
     \begin{subfigure}[b]{0.45\textwidth}
        \centering
        \includegraphics[width=5cm]{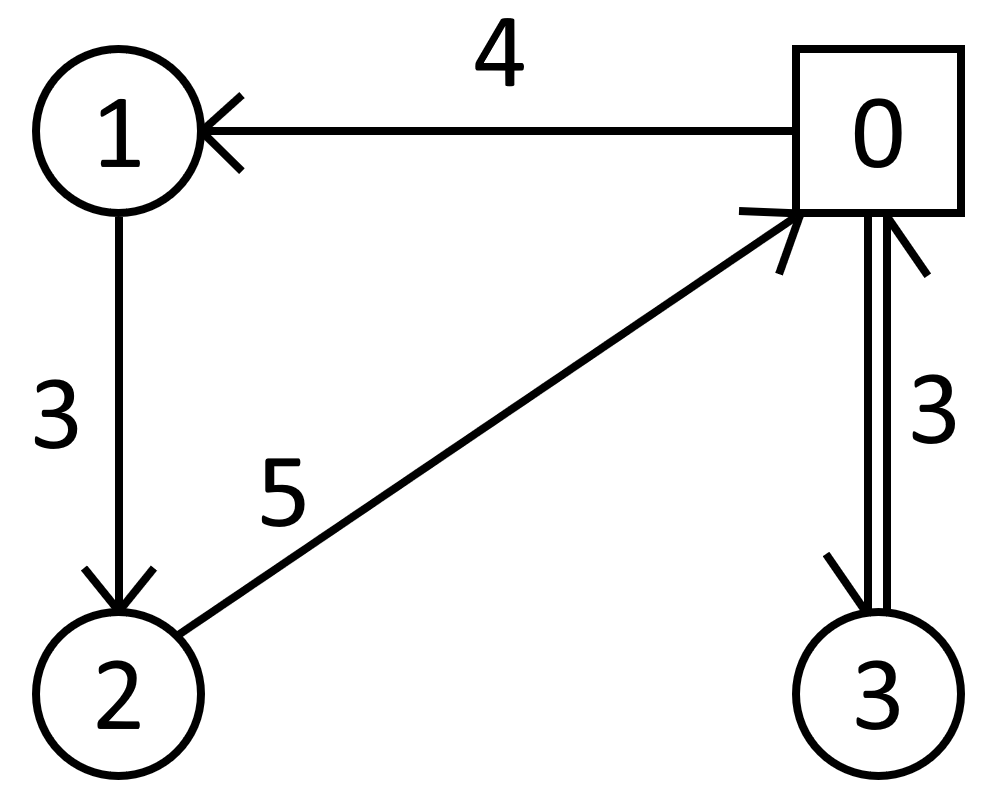}
        \caption{\footnotesize Routing and clustering of the example solution.}
        \label{illusRoutes}
     \end{subfigure} 
     \caption{Locations (a) and the associated routing (b) for the illustrative SCIRP example.}
\end{figure}
The tactical-level decision is displayed in Figure \ref{illusRoutes}. It shows two clusters being selected: A cluster comprising customers $1$ and $2$, and a cluster consisting of customer 3 solely. The transportation delivery schedules and associated costs are displayed in Table \ref{illuSolnData}. 

\begin{table}[h]
     \caption{The transportation delivery schedule and the costs of the illustrative solution.}
     \label{illuSolnData}
    \begin{subtable}[h]{0.45\textwidth}
        \centering
        \caption{\footnotesize The transportation delivery schedule of clusters.}
        \label{illuSchedules}
        \begin{tabular}{l|ccccccc} \toprule
        Cluster & Mo & Tu & We & Th & Fr & Sa & Su \\ \midrule
        $\{1,2\}$ & $\checkmark$ &  &  & $\checkmark$ &  & $\checkmark$ &  \\
        $\{3\}$ & & $\checkmark$ & $\checkmark$ & & $\checkmark$ & & $\checkmark$ \\ \bottomrule
        \end{tabular}
    \end{subtable}
    \hfill
    \begin{subtable}[h]{0.45\textwidth}
        \centering
        \caption{\footnotesize The cycle costs of clusters.}
        \label{illuCosts}
        \begin{tabular}{l|rrr} \toprule
        Cluster & Transport & Holding & Emergency \\ \midrule
        \{1,2\} & 1020      & 863    & 72.2     \\
        \{3\}   & 880      & 885    & 77.4 \\\bottomrule  
        \end{tabular}
     \end{subtable}
\end{table}

It is found that the transportation cost of cluster $1$ is 1020. With selected parameters $W = 100$ and $w = 20$ and the tour length of $4+3+5=12$, one replenishment of cluster $1$ costs $100 + 20 \times 12 = 340$. The schedule implies having $3$ visits weekly on cluster $1$, totaling a cyclic transportation cost of $340 \times 3 = 1020$. Secondly, the base-stock levels for cluster with $\alpha = 0.95$ are selected as $\{371,258,258\}$ for customer $1$ and $\{892,616,616\}$ for customer $2$, respectively on Monday, Thursday, and Saturday replenishments. These base-stocks does not exceed the customer inventory capacity $U=1200$. Expected positive amount of average inventories of customers $1$ and $2$ are calculated with these base-stocks as in Table \ref{illusHold}. On average, with $h = 0.2$, a holding cost of $0.2 \times 185.35 = 37.07$ and $0.2 \times 431.42 = 86.22$ is paid daily for customers $1$ and $2$, respectively. This results in an expected holding cycle cost of $(86.22 + 37.07) \times 7 = 863$ for the cluster. Lastly, to calculate emergency shipment amounts, we provide the mean and the standard deviation of normal distributions for emergency shipment quantities, and corresponding expected positive part of these distributions in Table \ref{illusEmerg} for $Q = 1200$. Assuming $e = 25$, the expected emergency shipment cost would be $25 \times 2.89 = 72.2$ on cycle.

\begin{table}[h]
     \caption{The transportation delivery schedule and the costs of the illustrative solution.}
     \label{illusResultData}
    \begin{subtable}[h]{0.45\textwidth}
        \centering
        \caption{\footnotesize Expected average inventory of cluster $1$.}
        \label{illusHold}
        \begin{tabular}{r|rrrrrrr|r} \toprule
        Cust. & Mo & Tu & We & Th & Fr & Sa & Su & Avg\\ \midrule
        $1$ & 321 & 221 & 122 & 208 & 109 & 208 & 109 & 185\\
        $2$ & 767 & 517 & 268 & 491 & 242 & 491 & 242 & 431\\\bottomrule
        \end{tabular}
    \end{subtable}
    \hfill
    \begin{subtable}[h]{0.45\textwidth}
        \centering
        \caption{\footnotesize Distribution of emergency shipments.}
        \label{illusEmerg}
        \begin{tabular}{l|rrr} \toprule
        Day & Mean & Std. Dev. & Exp. Pos. \\ \midrule
        Mo & -110.79  & 79.06    & 2.89     \\
        Th   & -539.21      & 96.82    & 0.00 \\
        Sa   & -500.00      & 79.06    & 0.00 \\
        \bottomrule  
        \end{tabular}
     \end{subtable}
\end{table}

Upon the transportation delivery schedule, we need to base our dynamic purchasing decisions. We discretize the continuous inventory levels in steps of $5$ units. We set $K_1 = 1000$, $K_2 = 2500$, $b_1 = 10$, and $b_2 = 2$. The resulting expected cyclic buying and selling costs, obtained by value iteration, is found to be $782$. We observe that the optimal purchase policy is similar to a period-dependent $(s,S)$ policy, i.e., if the inventory drops below a corresponding $s$, the policy buys up to $S$. The $(s,S)$ values for each day of the week are provided in Table \ref{illusS}.

\begin{table}[h]
\centering
\caption{The optimal $(s,S)$ policies of the example solution.}
\label{illusS}
\begin{tabular}{r|rrrrrrr} \toprule
 & Mo & Tu & We & Th & Fr & Sa & Su \\ \midrule
s & 0 & 340 & 30 & 315 & 70 & 350 & 330 \\ 
S & 385 & 795 & 605 & 805 & 670 & 820 & 640  \\\bottomrule
\end{tabular}
\end{table}

Unlike in a classic stationary $(s,S)$ policy, the producer's inventory may exceed $S$ for particular periods, for example, when the supply is more than the demand outflow on the next day of making a purchase. Moreover, this may occur consecutively for multiple periods and accumulate the producer's inventory up to its capacity. In the ranges of $[s,C]$ for each day of the week, the policy implies no buy or sell for the producer in the illustrative example. If the inventory is above $C$, the optimal policy is only to sell the part that is in excess of the capacity. 

This solution is efficient where both decision stages are considered jointly. 
In order to illustrate this efficiency, and the idea behind our approach in Section \ref{sect:soln}, consider an alternative tactical-level solution where the only change occurs on the transportation delivery schedule of the clusters, given in Table \ref{illuSchedulesUpd}. In this solution, the transportation delivery schedule of customer $3$ is postponed for $5$ days. For example, the base-stock policy on Tuesday in the first solution is used as the base-stock policy on Sunday in the alternative solution. Thus, the expected cost of the tactical-level decision does not change by postponing the transportation delivery schedule. However, we observe that the expected cost of the dynamic purchasing decisions increases to $2025$, which is more than double the expected cost of the dynamic purchasing decision in the original solution. This results in a 27\% increase in the total expected cycle cost of the SCIRP.

\begin{table}[h]
     \centering
     \caption{The transportation delivery schedule of the alternative tactical-level solution.}
     \label{illuSchedulesUpd}
     \begin{tabular}{l|ccccccc} \toprule
        Cluster & Mo & Tu & We & Th & Fr & Sa & Su \\ \midrule
        $\{1,2\}$ & $\checkmark$ &  &  & $\checkmark$ &  & $\checkmark$ &  \\
        $\{3\}$ & $\checkmark$ &  & $\checkmark$ & & $\checkmark$ & & $\checkmark$ \\ \bottomrule
        \end{tabular}

\end{table}

\clearpage

\section{Pseudo-codes of the Solution Approach} \label{appendix2}

\begin{algorithm}[h]
\small
    \caption{Precompute the set $\mathscr{R}$.}\label{alg:setR}
    \hspace*{\algorithmicindent} \textbf{Input:} Demand distributions and instance parameters. \\
    \hspace*{\algorithmicindent} \textbf{Output:} The set $\mathscr{R}$ with corresponding parameters $c_r$, $\beta_r^i$, $\Delta_r^t$, and $\Lambda_r^t$.
    \begin{algorithmic}[1]
        \State Compute $F^{-1}_{in}(\alpha)$ for all $i \in \mathscr{N}$ and $n \in \mathscr{T}$.
        \State While computing $F^{-1}_{in}(\alpha)$, store $\Theta_i \coloneqq \max \{n \mid F^{-1}_{in}(\alpha) \leq U\}$ for all $i \in \mathscr{N}$.
        \State Recursively find all possible partitioning of customers, $\mathscr{K}$, in which a full transportation delivery schedule is possible without violating Constraint \eqref{cc2}.
        \State While computing $\mathscr{K}$, store $\mathit{UB}^k$ as the longest time in between two consecutive deliveries without violating Constraint \eqref{cc2} for partition $k \in \mathscr{K}$.
        \ForEach {partition $k \in \mathscr{K}$}
        \ForEach {$2^{T}-1$ transportation delivery schedules}
        \If{the longest replenishment of the schedule  $ \leq \min \{ \mathit{UB}^k, \min_{i \in N^k} \Theta_i\}$}
        \State Add a cluster in  $\mathscr{R}$ by computing $c_r$, $\beta_r^i$, $\Delta_r^t$, and $\Lambda_r^t$ for all $i \in \mathscr{N}$ and $t \in \mathscr{T}$.
        \EndIf
        \EndFor
        \EndFor
    \end{algorithmic}
\end{algorithm}

\begin{algorithm}[h]
\small
    \caption{Line Search.}\label{alg:line}
    \hspace*{\algorithmicindent} \textbf{Input:} Initial values of $(\eta_1,\eta_2) = (\epsilon,\epsilon)$, increment value of $\zeta$, an upper bound of $UB$, and the set $\mathscr{R}$ with corresponding parameters $c_r$, $\beta_r^i$, $\Delta_r^t$, and $\Lambda_r^t$. \\
    \hspace*{\algorithmicindent} \textbf{Output:} The minimum found cost of $z = \sum_{r \in \mathscr{R}} c_r x_r  + C(\pi(\mathbf{x}))$.
    \begin{algorithmic}[1]
        \State Initialize $z = \infty$, $\psi = 0$, and $i = 0$.
        \While {$i \ne 2$ and $(\eta_1,\eta_2) \ne (UB,UB)$}
        \State Solve the modified set-partitioning model with $c_r$, $\beta_r^i$, $\Delta_r^t$, $\Lambda_r^t$, and $(\eta_1,\eta_2)$ to record $x'_r$ and $\sum_{r \in \mathscr{R}} c_r x'_r$.
        \State Solve MDP with the corresponding tactical routing solution, $x'_r$, to derive the optimal policies.
        \State Simulate the optimal MDP policy to derive $C(\pi(\mathbf{x'}))$ and record $z' = \sum_{r \in \mathscr{R}} c_r x'_r + C(\pi(\mathbf{x'}))$.
        \If{$z' \leq z$}
        \State Update $z \coloneqq z'$.
        \If{$\max(\eta_1,\eta_2) < UB$}
        \State Update $\eta_1 \coloneqq \eta_1 + \psi \zeta$, and $\eta_2 \coloneqq \eta_2 + (1 - \psi) \zeta$.
        \Else
        \State Update $\psi \coloneqq 1-\psi$ and $i \coloneqq i + 1$.
        \EndIf
        \Else
        \State Update $\eta_1 \coloneqq \eta_1 - \psi \zeta$, and $\eta_2 \coloneqq \eta_2 - (1 - \psi) \zeta$.
        \State Update $\psi \coloneqq 1-\psi$ and $i \coloneqq i + 1$.
        \EndIf
        \EndWhile
    \end{algorithmic}
\end{algorithm}

\clearpage

\section{Case Study Information} \label{appendix1}
 
\begin{table}[h]
\caption{Mean and standard deviations of daily demand distributions (in kg) of the case study.} \label{AllDist}
\begin{tabular}{ll||rr||rr||rr} \toprule
     &              & Emmen &         & Eems1 &         & Eems2 &         \\
Node & Location     & Mean & Std Dev & Mean & Std Dev & Mean & Std Dev \\ \midrule 
0    & Emmen        & 500  & 250     & -    & -       & -    & -       \\
0    & Eemshaven    & -    & -       & 2195 & 658.5   & 4670 & 934     \\\midrule
1    & Groningen    & 100  & 20      & 350  & 52.5    & 600  & 60      \\
2    & Delfzijl     & 35   & 14      & 100  & 25      & 150  & 30      \\
3    & Leeuwarden   & 100  & 20      & 200  & 40      & 400  & 60      \\
4    & Emmen        & 65   & 16      & 250  & 50      & 450  & 67.5    \\
5    & Pesse        & 50   & 15      & 120  & 30      & 250  & 37.5    \\
6    & Meendenertol & 30   & 12      & 75   & 22.5    & 120  & 24      \\
7    & Panjerd      & 30   & 12      & 100  & 30      & 150  & 30      \\
8    & Haerst       & 30   & 12      & 100  & 30      & 150  & 30      \\
9    & Paardeweide  & 30   & 12      & 100  & 30      & 150  & 30      \\
10   & Oude Riet    & 30   & 12      & 100  & 30      & 150  & 30      \\
11   & Roode Til    & -    & -       & 100  & 30      & 150  & 30      \\
12   & Veenborg     & -    & -       & 100  & 30      & 150  & 30      \\
13   & Smalhorst    & -    & -       & 100  & 30      & 150  & 30      \\
14   & Zeijerveen   & -    & -       & 100  & 30      & 150  & 30      \\
15   & Mandelân     & -    & -       & 100  & 30      & 150  & 30      \\
16   & De Horne     & -    & -       & 100  & 30      & 150  & 30      \\
17   & Stienkamp    & -    & -       & 100  & 30      & 150  & 30      \\
18   & Bloksloot    & -    & -       & -    & -       & 150  & 30      \\
19   & De Krellen   & -    & -       & -    & -       & 150  & 30      \\
20   & Mienscheer   & -    & -       & -    & -       & 150  & 30      \\
21   & Glimmermade  & -    & -       & -    & -       & 150  & 30      \\
22   & Dekkersland  & -    & -       & -    & -       & 150  & 30      \\
23   & Bareveld     & -    & -       & -    & -       & 150  & 30      \\
24   & Dikke Linde  & -    & -       & -    & -       & 150  & 30      \\ \bottomrule
\end{tabular}
\end{table}

\end{APPENDICES}

\end{document}